\documentclass[11pt]{amsart}

\usepackage{amssymb, amsmath}
\allowdisplaybreaks
\usepackage[all]{xy}
\usepackage[inline]{enumitem}
\usepackage{hyphenat}
\usepackage[colorlinks,linkcolor=blue,citecolor=blue,urlcolor=teal]{hyperref}
\usepackage{xcolor}
\usepackage{mathtools}
\usepackage{tikz-cd}
\usepackage[labelformat=empty]{caption}
\usepackage{cleveref} 

\usepackage[ left=3cm, right=3cm, top = 2.8cm, bottom = 2.8cm ]{geometry}

\newcommand{\eq}[2]{\begin{equation}\label{#1}#2 \end{equation}}

\newcommand{\xr}[1] {\xrightarrow{#1}}

\newcommand{\ra}{\rightarrow}

\newcommand{\N}{\mathbb{N}}

\newcommand{\Z}{\mathbb{Z}}

\newcommand{\sB}{\mathfrak{B}}
\newcommand{\sE}{\mathfrak{E}}

\newcommand{\sD}{\mathfrak{D}}

\newcommand{\cO}{\mathcal{O}}

\newcommand{\cZ}{\mathcal{Z}}

\newcommand{\ssm}{\mathfrak{m}}

\newcommand{\Ker}{\operatorname{Ker}}

\renewcommand{\Im}{\operatorname{Im}}

\newcommand{\Frac}{\operatorname{Frac}}

\newcommand{\xra}{\xrightarrow}
\renewcommand{\lim}{\operatornamewithlimits{\varprojlim}}

\newcommand{\Spec}{\operatorname{Spec}}

\newcommand{\Div}{\operatorname{Div}}

\newcommand{\red}{{\operatorname{red}}}

\newcommand{\dlog}{\operatorname{dlog}}

\newcommand{\hra}{\hookrightarrow}
\renewcommand{\hat}{\widehat}

\newcommand{\Zar}{{\operatorname{Zar}}}
\newcommand{\Nis}{{\operatorname{Nis}}}
\newcommand{\et}{{\operatorname{\acute{e}t}}}

\renewcommand{\phi}{{\varphi}}

\newcommand{\ul}{\underline}
\renewcommand{\tilde}{\widetilde}

\renewcommand{\bar}{\overline}

\theoremstyle{plain}

\newtheorem{prop}{Proposition}[section]
\newtheorem{lem}[prop]{Lemma}

\newtheorem{cor}[prop]{Corollary}
\newtheorem{thm}[prop]{Theorem}
\newtheorem{claim}{Claim}[prop]
\newtheorem{Observation}{Observation}[prop]

\theoremstyle{definition}
\newtheorem{defn}[prop]{Definition}

\newtheorem{example}[prop]{Example}

\newtheorem{para}[prop]{}

\newtheorem{remark}[prop]{Remark}

\newcommand{\beq}{\begin{equation}}
\newcommand{\eeq}{\end{equation}}

\numberwithin{equation}{prop}

\title{Structure theorem for log de Rham-Witt sheaves with vanishing}
\author{Fei Ren}

\address{Bergische Universit\"at Wuppertal\\ Gau\ss str. 20, D-42119 Wuppertal, Germany}
\email{renfei@uni-wuppertal.de}

\thanks{ }
\begin{document}
\begin{abstract}
The $\Z/p^n$ coefficient motivic cohomology with vanishing along an effective divisor $D$ is characterized by log de Rham-Witt sheaves with vanishing along the same divisor.
In this paper, we prove an elegant structure theorem for log de Rham-Witt sheaves with vanishing along $D$ defined in \cite{RamFil1}, answering a question of Shuji Saito during the Mainz conference and a question of Yigeng Zhao during a short visit of the author last summer.
Our structural result for the log forms also lays the foundation for the study of Milnor $K$-theory with vanishing along $D$ in the paper to come. 
\end{abstract}
\maketitle

\tableofcontents

\section{Introduction}
Let $k$ be a perfect field of positive characteristic $p$.
Let $X$ be a smooth $k$-scheme (not necessarily proper) of pure dimension $d$ and let $D$ be
an effective Cartier divisor such that the underlying reduced divisor $D_{\red}$ has simple normal
crossings.
Let $U$ be the complement of $D$ in $X$ and $j:U\hra X$ be the open immersion.
In \cite{RamFil1}, we proved a duality theorem for de Rham-Witt sheaves with certain pole- or zero-restrictions along $D$ in the Nisnevich topology. 
The \textit{de Rham-Witt sheaf with zeros along $D$}, defined by
$$W_n\Omega^q_{(X,-D)}:=\Ker\left(W_n\Omega^q_X\ra\bigoplus_\alpha W_n\Omega^q_{n_{\alpha} \sD_\alpha}\right)
 \quad$$
(where $D=\sum_\alpha n_{\alpha} \sD_\alpha$ with all $\sD_{\alpha}$
\footnote{The notation $\sD_\alpha$ here is different from the notation $D_i$ in \Cref{thm0}, so we use different fonts as well as different lower indices to emphasize this. $\sD_\alpha$ here are integral, while $D_i$ that appears in a $p$-divisibility decomposition are in general neither irreducible nor reduced. 
 }
being the irreducible components of $D$ and $D_{\red}=\sum \sD_\alpha$ is  strict normal crossing)
is dual to 
the the corresponding sheaf with poles defined via the theory of sheaves with modulus \cite{KMSYI, KMSYII, KMSYIII}.
This duality extends further to a duality theorem of Milne-Kato type.
Namely, the \textit{$\Z/p^n$-coefficient motivic complex with zeros along $D$}, defined by
\[\Z/p^n(q)_{(X,-D)}:=\left(W_n\Omega^q_{(X,-D)}\xr{C^{-1}-1} 
\frac{W_n\Omega^q_{(X,-D)}}{dV^{n-1}\Omega^{q-1}_X\cap W_n\Omega^q_{(X,-D)}}\right)[-q],\]
is dual to the corresponding motivic complex with poles in the \'etale topology.
See \cite[Theorem 11.15]{RamFil1} for the precise statements.
Moreover, the complex $\Z/p^n(q)_{(X,-D)}$ is concentrated in degree $q$ in the \'etale topology by \cite[Lemma 11.6]{RamFil1}. This only nonzero cohomology sheaf is denoted by $W_n\Omega^q_{(X,-D),\log}$, and it is the intersection of the usual log de Rham-Witt sheaf $W_n\Omega^q_{X,\log}$ with 
$W_n\Omega^q_{(X,-D)}$.
We will refer to it as the \textit{log de Rham-Witt sheaf with vanishing along $D$} or the \textit{log de Rham-Witt sheaf with zeros along $D$}.
Hence to study the $\Z/p^n$-motivic complex with zeros along $D$, it is essential to understand the structure of $W_n\Omega^q_{(X,-D),\log}$.
In this work, we prove the following fundamental structural result:
\begin{thm}\label{thm0}
Let $n\ge 1$ be an integer.
Let $$D=D_0+pD_1+\cdots +p^{n}D_{n}$$ be the $p$-divisibility decomposition of $D$ of length $n$ (see \S\ref{para:div-not} for the precise definition).
For each $i$, set
$$\ul D_i=D_0+pD_1+\cdots + p^iD_i.$$
Let
$$j_i:U_i=X\setminus \ul D_i\hra X,\qquad
j:U=X\setminus D\hra X$$ be the open immersions.
Let 
$$\dlog:j_*(\underbrace{\cO_{U}^\times \otimes \cdots \otimes \cO_{U}^\times)}_{q\text{-times}}
\ra j_*W_{n}\Omega^q_U,\quad 
a_1\otimes \cdots\otimes a_q\mapsto \frac{d[a_1]\cdots d[a_q]}{[a_1\cdots a_q]}$$
be the $\dlog$ map, where $[-]:=[-]_{n}$ is the Techm\"uller lift.
Then
$$W_n\Omega^q_{(X,-D),\log}=\sum_{i=0}^{n-1} 
 p^i \dlog\left(
 \Ker(\cO_X^\times\ra \cO_{ \lceil\frac{D}{p^i}\rceil}^ \times)\otimes
\underbrace{ j_{i*}\cO_{U_i}^\times \otimes \cdots \otimes j_{i*}\cO_{U_i}^\times}_{(q-1)\text{-times}} \right)$$
as subsheaves of $W_n\Omega^q_X$ on $X_\Nis$.
\end{thm}
We remind the reader that neither the $p$-divisibility decomposition given in the statement is  assumed to be of maximal length, nor any of $D_0,\cdots,D_n$ is assumed to be nonzero.
(In other words, the shape of $D$ is unrelated to the length $n$ appearing in $W_n\Omega^q$.)
To better illustrate our theorem \Cref{thm0}, we refer to \Cref{Ex1} and \Cref{Ex2} as special examples when $n=1$ and $D$ is of a particular form.

There are other definitions of log de Rham-Witt sheaves with zeros along $D$ earlier in the literature. 
For instance, \cite{Gupta-Krishna} defined a sheaf
$$\Ker(W_n\Omega^q_X\ra W_n\Omega^q_D)\cap W_n\Omega^q_{X,\log},$$
while \cite{JSZ} defined another sheaf
$$\Im\Big(\dlog: \Ker(\cO_X^\times\ra \cO_D^\times)
\otimes_{\Z} \underbrace{j_*\cO_U^\times\otimes_{\Z}\cdots \otimes_{\Z} j_*\cO_U^\times}_{(q-1)\text{-times}}
\ra W_n\Omega^q_X\Big).$$
As $n$ varies, all the three definitions form pro-systems, and these three pro-systems are isomorphic. 
It remained as a question how much they differ as sheaves. 
As a corollary of our structural theorem, we are able to give a complete answer to this question.
See \Cref{remark}.

This structural result for log forms with vanishing also laid the foundation for the study of Milnor $K$-theory in the paper to come.
To shed some light, we end this introduction session by providing our definition of the Milnor $K$-sheaf with vanishing along $D$. Suppose $D=D_0+pD_1+\cdots+ p^LD_L$ is the $p$-divisibility decomposition of $D$ of maximal length (see \S\ref{para:div-not} for the precise definition), then we define
$$K^M_{q,(X,D)}:=
\sum_{i=0}^{L} 
\Im\left(
 p^i  \Ker(\cO_X^\times\ra \cO_{ \lceil\frac{D}{p^i}\rceil}^ \times)
\otimes_\Z 
j_{i*} K^M_{q-1,U_i}
\ra 
j_*K^M_{q,U}
\right)$$
as an abelian sheaf on $X_\Nis$.
With this definition, we arrive at a relative version of the Bloch-Gabber-Kato theorem.
\begin{cor}
The $\dlog$ map induces an isomorphism of sheaves
$$\dlog: K^M_{q,(X,D)}/p^n K^M_{q,X}\cap  K^M_{q,(X,D)}
\xra{\simeq}
W_n\Omega^q_{(X,-D),\log}$$
on $X_\Nis$.
\end{cor}

\subsection{Notations and conventions}
\label{para:div-not}
Throughout the article, $k$ denotes a perfect field of positive characteristic $p$. 
$X$ is a smooth separated $k$-scheme of finite type.
For an effective Cartier divisor $E$ on $X$ with $E_\red$ being a simple normal crossing divisor (SNCD), we denote
$\Omega^q_X(\log E):=\Omega^q_{X}(\log E_\red).$

If $E$ is any divisor on $X$ and $m\in \Z\setminus\{0\}$, then we write
\[m \mid E :\Longleftrightarrow  m \text{ divides the multiplicity of every irreducible component of } E\]
and 
\[m\nmid E :\Longleftrightarrow m \text{ does not divide the multiplicity of any irreducible component of } E.\]
Let $n\ge 0$ be an integer.
We say
\eq{eqpdd}{E= E_0 + pE_{1}+\ldots + p^{n}E_{n}
}
is a {\em $p$-divisibility decomposition  of $E$ (of length $n$)}, if
\begin{itemize} 
\item 
$p^{i}\nmid p^{i-1}E_{i-1}$, for $i=1,\ldots, n$, 
and 
\item 
$\sum_{i=0}^n E_{i,\red}$ is a reduced divisor.
\end{itemize}
Note that 
\begin{enumerate}
\item 
$p\mid E_{n}$ is allowed,
\item 
$E_0,\cdots, E_s$ are allowed to be zero divisors, and
\item 
for any $n$, a $p$-divisibility decomposition of $E$ of length $n$ always exists and is unique. (Hence we are allowed to use the article ``the'' when the length is specified.)
\item 
If $\eqref{eqpdd}$ is the $p$-divisibility decomposition of length $n$, then
\eq{eqpddn-1}{E= E_0 + p E_{1}+\ldots +p^{n-2}E_{n-2}+ p^{n-1}(E_{n-1}+pE_{n})
}
is a $p$-divisibility decomposition of length $n-1$.
\end{enumerate}
When $p\nmid E_n$, we say that \eqref{eqpdd} is the $p$-divisibility decomposition of $E$ \textit{of maximal length}.
For any given $n$, the $p$-divisibility decomposition of length $n$ in this sense is a $p$-divisibility decomposition with respect to the sequence $1<2<\cdots<n$ as defined in \cite{RamFil1}.
In this article, we will be only using $p$-divisibility of this shape, so it suffices to use the above definition.
If we have
$E=\sum_i n_i \sE_i$ with $\sE_i$ the irreducible components of $E$, $\sE_i\neq\sE_j$, for $i\neq j$,
then we set
\[\lceil E/m\rceil:= \sum_i \lceil n_i/m\rceil \sE_i \quad \text{and} \quad \lfloor E/m\rfloor:= \sum_i \lfloor n_i/m\rfloor \sE_i,\]
where $\lceil - \rceil$ (resp. $\lfloor -\rfloor$) denotes the  round-up (resp. round-down).

Finally we specify a few shorthands. By $p\mid r\ge 1$ we mean $p\mid r$ and $r\ge 1$. By $p\nmid r\ge 1$ we mean $p\nmid r$ and $r\ge 1$.

\subsection*{Acknowledgement}
The author thanks Kay R\"ulling for various discussions and for checking an earlier draft of this article. 
Gratitude also goes to Shuji Saito and Yigeng Zhao for raising questions that steered the exploration in this direction.

\section{The structure of log de Rham-Witt differentials}
Let $D_{\red}$ be a SNCD. 

\begin{para}
Assume $X$ is of pure dimension $d$.
Recall from \cite{RamFil1} that the \textit{de Rham-Witt sheaf with zeros along $D$} is defined by
$$W_n\Omega^q_{(X,-D)}:=\Ker\left(W_n\Omega^q_X\ra\bigoplus_\alpha W_n\Omega^q_{n_{\alpha} \sD_\alpha}\right),
$$
where $D=\sum_\alpha n_{\alpha} \sD_\alpha$ is the decomposition of $D$ into its irreducible components, namely all $\sD_{\alpha}$
being the irreducible components of $D$ and $D_{\red}=\sum \sD_\alpha$ is strict normal crossing.

\begin{enumerate}[label=(\arabic*)]
\item 
$W_n\Omega^q_{(X,-D)}$ is Cohen-Macaulay (see \Cref{lemAppHartog:Ex}\eqref{lemAppHartog:Ex3}).
\item 
All the maps $R,\ul p,F,V,d$ of $W_n\Omega^\bullet_X$ restricts to well-defined endomorphisms of $W_n\Omega^\bullet_{(X,-D)}$. 
More precisely, 
$$R: W_n\Omega^q_{(X,-D)}\ra W_{n-1}\Omega^q_{(X,-D)},\quad
F:R: W_n\Omega^q_{(X,-D)}\ra W_{n-1}\Omega^q_{(X,-D)}$$
$$\ul p: W_n\Omega^q_{(X,-D)}\ra W_{n+1}\Omega^q_{(X,-D)},\quad
V:W_n\Omega^q_{(X,-D)}\ra W_{n+1}\Omega^q_{(X,-D)}$$
$$d: W_n\Omega^q_{(X,-D)}\ra  W_n\Omega^{q+1}_{(X,-D)}$$
are well-defined.
\item 
The inverse Cartier operator 
$$C^{-1}: W_n\Omega^q_{X} \ra
\frac{W_n\Omega^q_{X}}{dV^{n-1}\Omega^{q-1}_X}$$
restricts to a well-defined map
$$C^{-1}:W_n\Omega^q_{(X,-D)}\ra 
\frac{W_n\Omega^q_{(X,-D)}}{dV^{n-1}\Omega^{q-1}_X\cap W_n\Omega^q_{(X,-D)}}.$$

\end{enumerate}

\end{para}

\begin{defn}[{cf. \cite[Lemma 11.6]{RamFil1}}]
Define $W_n\Omega^q_{(X,-D),\log}$ to be the abelian subsheaf of $W_n\Omega^q_{X,\log}$ on $X_\Nis$ such that the following sequence is exact
$$0\ra
W_n\Omega^q_{(X,-D),\log}\ra
W_n\Omega^q_{(X,-D)}\xra{C^{-1}-1}
\frac{W_n\Omega^q_{(X,-D)}}{dV^{n-1}\Omega^{q-1}_X\cap W_n\Omega^q_{(X,-D)}}.$$
We write $\Omega^q_{(X,-D),\log}:=W_1\Omega^q_{(X,-D),\log}$.
\end{defn}
\begin{remark}
As 
$$W_n\Omega^q_{X,\log}=
\Ker\left( W_n\Omega^q_{X}\xra{C^{-1}-1}
\frac{W_n\Omega^q_{X}}{dV^{n-1}\Omega^{q-1}_X\cap W_n\Omega^q_{X}}
\right)$$
by \cite[Lemma 4.1.5]{Kato-dualityI},
we see
$W_n\Omega^q_{(X,-D),\log}=W_n\Omega^q_{(X,-D)}\cap W_n\Omega^q_{X,\log}$.
\end{remark}

\begin{thm}
\label{thm1}
Let $A, B\ge 0$ be Cartier divisors on $X$ with $A$ and $(A+B)_\red$ being SNCDs. 
Write $$B=B'+B''$$ 
where 
$B'_\red\le A$, $B''$ and $A$ share no irreducible components, and $B'_\red+B''_\red$ being a SNCD. 
Let $$B''=B''_0+pB''_1$$
be the $p$-divisibility decomposition of $B''$ of length 1. 
Denote by $j:U=X\setminus (A+B)\hra X$ the open immersion.
For $q\ge 1$, let
$$\Omega^q_X(\log A)(-B)_{\log}:=
j_*\Omega^q_{U,\log}\cap \Omega^q_X(\log A)(-B).
$$
Let
$$\dlog:j_*(\underbrace{\cO_{U}^\times \otimes \cdots \otimes \cO_{U}^\times}_{q\text{-times}})
\ra j_*\Omega^q_U,\quad 
a_1\otimes \cdots\otimes a_q\mapsto \frac{da_1\wedge \cdots\wedge da_q}{a_1\cdots a_q}$$
be the $\dlog$ map.
Then 
we have the following identity of Nisnevich subsheaves of $\Omega^q_X(\log (A+B))$
$$\Omega^q_X(\log A)(-B)_{\log}=
\dlog
\left(
\Ker(\cO_X^\times \ra \cO_{\tilde B}^\times)\otimes 
\underbrace{ j_{0*}\cO_{V}^\times\otimes\cdots \otimes  j_{0*}\cO_{V}^\times}_{(q-1)\text{ times}}
\right).$$
where 
$$\tilde B=B+B''_{0,\red},$$
and
$j_0:V:=X\setminus (A+B''_0)\hra X$ is the open immersion.
\end{thm}
\begin{cor}\label{corn=1}
Let $D=D_0+pD_1$ is the $p$-divisibility decomposition of $D$ of length 1.
Recall that $\Omega^q_{(X,-D)}=\Omega^q_X(\log D_0)(-D)$ by \cite[Lemma 8.4]{RamFil1}.
Therefore \Cref{thm1} directly gives
$$\Omega^q_{(X,-D),\log}=
\dlog\left(
\Ker(\cO_X^\times \ra \cO_{D}^\times)\otimes 
\underbrace{ j_{0*}\cO_{U_0}^\times\otimes\cdots \otimes  j_{0*}\cO_{U_0}^\times}_{(q-1)\text{ times}}
\right)$$
where $j_0:U_0:=X\setminus D_0\hra X$ is the open immersion.
\end{cor}
Before discussing the proof of \Cref{thm1}, we give two examples of \Cref{corn=1} when the $p$-divisibility decomposition of $D$ is of a particular shape.
\begin{example}\label{Ex1}
Suppose $D=D_0$ is the $p$-divisibility decomposition of $D$  of length 1, 
i.e., multiplicities of all irreducible components of $D$ are coprime to $p$.
In this case, we have $\Omega^q_{(X,-D)}=\Omega^q_X(\log D)(-D)$, and it follows from \Cref{corn=1} that
$$\Omega^q_{(X,-D),\log}=
\dlog\left(
\Ker(\cO_X^\times \ra \cO_{D}^\times)\otimes 
\underbrace{ j_{*}\cO_{U}^\times\otimes\cdots \otimes  j_{*}\cO_{U}^\times}_{(q-1)\text{ times}}
\right),$$
where $j:U:=X\setminus D\hra X$ is the open immersion.
\end{example}
\begin{example}\label{Ex2}
Suppose $D=pD_1$ is the $p$-divisibility decomposition of $D$ of length 1, i.e., multiplicities of all irreducible components of $D$ are divisible by $p$.
In this case, we have $\Omega^q_{(X,-D)}=\Omega^q_X(-D)$, and it follows from \Cref{corn=1} that
$$\Omega^q_{(X,-D),\log}=
\dlog\left(
\Ker(\cO_X^\times \ra \cO_{D}^\times)\otimes 
\underbrace{ \cO_{X}^\times\otimes\cdots \otimes  \cO_{X}^\times}_{(q-1)\text{ times}}
\right)$$
\end{example}

\begin{proof}[{Proof of \Cref{thm1}}]
Since $X$ is of finite type over $k$, it suffices to check that the two subsheaves have the same Nisnevich stalk at each closed point of $X$.
Let $R$ be the henselization of a local ring at a closed point of $X$.
\Cref{thm1} now follows from \Cref{MainProp}.
\end{proof}

\begin{prop}
\label{MainProp}
Let $(R,\ssm,k')$ be the henselization of a closed point $x$ in $X$, where $\ssm$ is the maximal ideal of $R$, and $k'$ is the residue field.
Let $T_1,\cdots , T_d$ be a system of regular parameters of $R$ such that 
$$A=\Div(T_1\cdots T_e), \quad
B=\Div(T_1^{r_1}\cdots   T_g^{r_g})$$
where $e,f,g$ are integers such that $0\le e\le f \le g\le d$  and $r_1,\cdots, r_g$ are nonnegative integers such that at least one $r_i\ge 1$, and
$$\text{$p\nmid r_j\ge 1$ for all $j\in [e+1,f]$,
\quad 
$p\mid r_j\ge 1$ for all $j\in [f+1,g]$.}$$
Put 
$r_{g+1}=\cdots =r_{d}=0$ and set
$$\ul r:=(r_1,\cdots, r_d)\in \N^d.$$
Set
$$G_{R}^{\ul r}:=
\Omega^q_R(\log A)(-B),\quad
G_{R,\log}^{\ul r}:=
\Omega^q_R(\log A)(-B)
\cap
\Omega^q_{\Frac R,\log}.$$
Define
$$\tilde r_i=
\begin{cases}
    r_i+1& i\in[e+1,f];
    \\
    r_i& i\notin [e+1,f].
\end{cases}$$
Set $\tilde{\underline{r}}:=(\tilde r_1,\cdots,\tilde r_d)$.
Then 
\begin{align}
\label{MainPropeq0}
G_{R,\log}^{\ul r}=
\Big\{
\sum_j^{\text{finite}}\dlog x_{1,j}\wedge &\dlog x_{2,j}\wedge\cdots \wedge\dlog x_{q,j}
\Big|
\\
&x_{1,j}\in 
1+(T_1^{\tilde r_1}\cdots  T_g^{\tilde r_g})\cdot R,\,
x_{2,j},\cdots,x_{q,j} \in R[\frac{1}{T_1\cdots T_f}]^\times
\Big\}.
\nonumber
\end{align}
\end{prop}
\Cref{MainProp} is the main technical part of this work. Its proof is long and tedious and we postpone it to the next session.

\begin{thm}
\label{thm2}
Let $n\ge 1$ be a given integer.
Let $$D=D_0+pD_1+\cdots +p^{n}D_{n}$$ be the $p$-divisibility decomposition of $D$ of length $n$.
For each $i$, set
$$\ul D_i=D_0+pD_1+\cdots p^iD_i.$$
Let
$$j_i:U_i=X\setminus \ul D_i\hra X,\qquad
j:U=X\setminus D\hra X$$ be the open immersions.
Let 
$$\dlog:j_*(\underbrace{\cO_{U}^\times \otimes \cdots \otimes \cO_{U}^\times)}_{q\text{-times}}
\ra j_*W_{n}\Omega^q_U,\quad 
a_1\otimes \cdots\otimes a_q\mapsto \frac{d[a_1]\cdots d[a_q]}{[a_1\cdots a_q]}$$
be the $\dlog$ map, where $[-]:=[-]_{n}$ is the Techm\"uller lift.
Then
\begin{enumerate}
\item 
We have the following identification of subsheaves of $W_n\Omega^q_X$ on $X_\Nis$:
$$W_n\Omega^q_{(X,-D),\log}=\sum_{i=0}^{n-1} 
  p^i \dlog\left(
 \Ker(\cO_X^\times\ra \cO_{ \lceil\frac{D}{p^i}\rceil}^ \times)\otimes
\underbrace{ j_{i*}\cO_{U_i}^\times \otimes \cdots \otimes j_{i*}\cO_{U_i}^\times}_{(q-1)\text{-times}} \right).$$
\item 
The restriction map
$$R:W_{n+1}\Omega^q_{(X,-D),\log}\ra W_n\Omega^q_{(X,-D),\log}$$
is surjective.
\end{enumerate}
\end{thm}

\begin{proof}
Once we have part (1), the surjectivity of $R$ clearly follows from the formula.
Hence our aim is to show the equality
$$W_n\Omega^q_{(X,-D)}\cap W_n\Omega^q_{X,\log}
=
\sum_{i=0}^{n-1} 
 p^i \dlog\left(
 \Ker(\cO_X^\times\ra \cO_{ \lceil\frac{D}{p^i}\rceil}^ \times)\otimes
\underbrace{ j_{i*}\cO_{U_i}^\times \otimes \cdots \otimes j_{i*}\cO_{U_i}^\times}_{(q-1)\text{-times}} \right).$$
When $n=1$, this is done in \Cref{corn=1}.
Now we assume $n\ge 2$ and denote the right hand side of the formula above by 
$M_n.$

We first show the direction $\supset$, that is
\eq{prop:Cinv-1n:eq1}{
\sum_{i=0}^{n-1} 
 p^i \dlog\left(
 \Ker(\cO_X^\times\ra \cO_{ \lceil\frac{D}{p^i}\rceil}^ \times)\otimes
\underbrace{ j_{i*}\cO_{U_i}^\times \otimes \cdots \otimes j_{i*}\cO_{U_i}^\times}_{(q-1)\text{-times}} \right)
\subset W_n\Omega^q_{(X,-D)}.
}
Since $W_n\Omega^q_{(X,-D)}$ is Cohen-Macaulay when restricted to $X_\Zar$ (see \Cref{lemAppHartog:Ex}\eqref{lemAppHartog:Ex3}), it suffices to check \eqref{prop:Cinv-1n:eq1} locally around all codimension $1$ points by \Cref{lemAppHartog}.
Let $R$ be the henselization of the local ring of a codimension $1$ point $x$ of $X$,  $\ssm$ be the maximal ideal of $R$, and $T$ be a uniformizer of $R$.
When $D$ does not pass through the point $x$, $\eqref{prop:Cinv-1n:eq1}$ trivially holds, as it is just the stalk at $x$ of the inclusion
$$\sum_{i=0}^{n-1} 
 p^i \dlog
\underbrace{ \cO_{X}^\times \otimes \cdots \otimes \cO_{X}^\times}_{q\text{-times}}
\subset W_n\Omega^q_{X}.$$
Now we suppose $D=\Div(T^r)$ in a neighborhood of $x$ with $r\ge 1$.
Suppose $r=p^jr'$ with $(p,r')=1$.
(Namely, $x$ is the generic point of an irreducible component of $D_j$ from the $p$-divisibility decomposition of $D$ of maximal length.)
The Nisnevich stalk of the left hand side of \eqref{prop:Cinv-1n:eq1} is
\eq{prop:Cinv-1n:eq3}{
\sum_{i=0}^{j-1} 
p^i \dlog\left(
(1+\ssm^{r/p^i})
\otimes R^\times \cdots \otimes R^\times\right)
+
\sum_{i=j}^{n-1} 
p^i \dlog \left(
(1+\ssm^{\lceil r/p^i\rceil})
\otimes R[\frac{1}{T}]^\times \cdots \otimes R[\frac{1}{T}]^\times
\right)
.
}
(When $j>n-1$, the second summand above is $0$.)
This is clearly contained in the Nisnevich stalk 
$$\Ker\left(W_n\Omega^q_R\ra  W_n\Omega^q_{R/\ssm^r}\right)$$
of the right hand side of \eqref{prop:Cinv-1n:eq1} at $x$.

We show the direction $\subset$.
The proof below is adapted from the proof of \cite[Theorem 2.3.1]{JSZ}.
We do induction on $n$.
Let $x\in W_n\Omega^q_{(X,-D)}\cap W_n\Omega^q_{X,\log}$ be a local section. 
Then $Rx\in W_{n-1}\Omega^q_{(X,-D)}\cap W_{n-1}\Omega^q_{X,\log}$, hence it lies in $M_{n-1}$ by induction hypothesis.
Consider the following commutative diagram
$$\xymatrix{
&&
M_n
\ar[r]^R\ar@{^(->}[d]&
M_{n-1}
\ar[r]\ar@{^(->}[d]&
0
\\
0\ar[r]&
\Omega^q_{X,\log}
\ar[r]^{\ul p^{n-1}}&
W_n\Omega^q_{X,\log}
\ar[r]^{R}&
W_{n-1}\Omega^q_{X,\log}
\ar[r]&
0.
}$$
Furthermore, 
$R:M_n\ra M_{n-1}$ is surjective by the very definition.
Let $y\in M_n$ such that $R(y)=R(x)$.
By the exactness of the second row (see \cite[Lemme 3]{CTSS}),
$x-y=\ul p^{n-1}(z)$ for some $z\in \Omega^q_{X,\log}$.
Note that $\ul p^{n-1}(z)=x-y\in W_n\Omega^q_{(X,-D)}$.
Hence by \cite[Lemma 8.11]{RamFil1}, we have 
$$z\in \Omega^q_{n-1}(-D',-pD_n):=\Omega^q_X(\log D')(-\lceil\frac{D}{p^{n-1}}\rceil),$$
where $D':=D-p^nD_n.$
In particular, $z\in \Omega^q_{X,\log}\cap \Omega^q_{n-1}(-D',-pD_n)$.
Note that
$$\lceil\frac{D}{p^{n-1}}\rceil= \lceil\frac{D'}{p^{n-1}}\rceil+pD_n
$$
is the $p$-divisibility decomposition of $\lceil\frac{D}{p^{n-1}}\rceil$ of length $1$, and $\Omega^q_X(\log D')=\Omega^q_X(\log \lceil\frac{D'}{p^{n-1}}\rceil)$ by definition (see \S\ref{para:div-not}).
Applying \Cref{thm1}, we have
$$z\in \dlog \left(
\Ker(\cO_X^\times\ra \cO_{\lceil\frac{D}{p^{n-1}}\rceil}^ \times)
\otimes \cO_{U_{n-1}}^\times \otimes\cdots \otimes \cO_{U_{n-1}}^\times
\right)$$
where $U_{n-1}=X\setminus D'=X\setminus\lceil\frac{D'}{p^{n-1}}\rceil$.
That is, $\ul p^{n-1}(z)\in M_n$.
\end{proof}

The following corollary is direct from the explicit formula for $W_n\Omega^q_{(X,-D),\log}$ given in \Cref{thm1}.
This is a filtered version of \cite[Lemme 3]{CTSS}.
Similar short exact sequences but for pro-systems were studied earlier for other filtrations, e.g. \cite[Proposition 5.9]{Gupta-Krishna} and \cite[Theorem 4.6]{Morrow}.
\begin{cor}
\label{Cor1}
The following sequence is exact
$$0\ra 
W_{n-1}\Omega^q_{(X,-D),\log}
\xra{\ul p}
W_n\Omega^q_{(X,-D),\log}
\xra{R^{n-1}}
\Omega^q_{(X,-D),\log}
\ra 0.
$$
\end{cor}

As an application for our structural result \Cref{thm1}, we answer a remaining question in \cite[11.5(4)]{RamFil1} and a question in \cite[Remark 11.7(1)]{RamFil1}.
In \cite{Gupta-Krishna} and \cite{JSZ} the following subsheaves of $W_n\Omega^q_{X}(\log D)$ are introduced:
$$W_n[{\rm GK}]^q_{(X,-D)}:=\Ker(W_n\Omega^q_X\ra W_n\Omega^q_D),\quad
W_n[{\rm GK}]^q_{(X,-D),\log}:=W_n[{\rm GK}]^q_{(X,-D)}\cap W_n\Omega^q_{X,\log};$$
$$\resizebox{\displaywidth}{!}{$
W_n[{\rm JSZ}]^q_{(X,-D)}:=W_n\cO_X(-D)\cdot W_n\Omega^q_X(\log D),\quad
W_n[{\rm JSZ}]^q_{(X,-D),\log}:=W_n[{\rm JSZ}]^q_{(X,-D)}\cap W_n\Omega^q_{X,\log}.$
}$$
Here $W_n\cO_X(-D)\subset W_n\cO_X$ refers to the ideal sheaf which is a free line bundle of $W_n\cO_X$-modules generated by the Teichm\"uller lifts of local sections of the effective Cartier divisor $D$.
\begin{cor}
\label{remark}
With notations as above, we have
\eq{Comp1}{W_n[{\rm GK}]^q_{(X,-D)}\subset 
W_n\Omega^q_{(X,-D)}\subset 
W_n[{\rm JSZ}]^q_{(X,-D)};}
\eq{Comp2}{W_n[{\rm GK}]^q_{(X,-D),\log}\subset 
W_n\Omega^q_{(X,-D),\log}\subset 
W_n[{\rm JSZ}]^q_{(X,-D),\log}.}
All the inclusions can be strict.
When $n=1$, we have the following more explicit expressions (we use the notations from \Cref{thm1}):
\eq{CompGK}{W_1[{\rm GK}]^q_{(X,-D)}
={\rm Im}(\cO_X(-D)\cdot \Omega^{q}_X+d\cO_X(-D)\cdot \Omega^{q-1}_X\hra \Omega^q_X),}
\eq{CompMein}{\Omega^q_{(X,-D)}=\Omega^q_X(\log D_0)(-D),\quad
\Omega^q_{(X,-D),\log}= 
\dlog\Ker(\cO_X^\times \ra \cO_{D}^\times)\wedge 
j_{0*}K^M_{q-1,U_0},}
\eq{CompJSZ}{W_1[{\rm JSZ}]^q_{(X,-D)}=\Omega^q_X(\log D)(-D),\quad 
W_1[{\rm JSZ}]^q_{(X,-D)}=
\dlog\Ker(\cO_X^\times \ra \cO_{D}^\times)\wedge 
\dlog j_{*}K^M_{q-1,U},}
\end{cor}
\begin{proof}
The equation \eqref{CompGK} follows from the fact that $W_n[{\rm GK}]^*_{(X,-D)}$ is the differential graded ideal of $W_n\Omega^*_{X}$ generated by $W_n\cO_X(-D)$, see e.g. \cite[Lemma 1.20]{Rue}.
The equations \eqref{CompMein} follow from \cite[Lemma 8.4]{RamFil1} and \Cref{thm1}.
The second equation of \eqref{CompJSZ} follows from \cite[Theorem 1.2.1]{JSZ} (or our \Cref{thm1}).
\begin{enumerate}
    \item 
Both inclusions in \eqref{Comp1} can be strict, this answers \cite[11.5(4)]{RamFil1}.
We illustrate this by the following explicit example when $n=1$. 
Let $X=\Spec A$,  $D=\Div(f^rg^sh^t)$
with $p\nmid rs$, $p\mid t$, and $D_{\red}=\Div(f)+\Div(g)+\Div(h)$ being an SNCD.
Then
\begin{align*}
\Omega^q_X(-D)\subsetneqq 
{\rm Im}(d\cO_X(-D)\cdot &\Omega^{q-1}_X\hra \Omega^q_X)
=
W_1[{\rm GK}]^q_{(X,-D)}
\subsetneqq \Omega^q_{(X,-D)}
\subsetneqq 
W_1[{\rm JSZ}]^q_{(X,-D)}
\end{align*}
where $M(-D)=M\otimes_{\cO_X} \cO_X(-D)$ for any $\cO_X$-module $M$.
Concretely, let $\omega\in\Omega^{q-1}_X$ be any nonzero regular differential, then
$$d(f^rg^sh^t)\wedge \omega\in 
{\rm Im}(d\cO_X(-D)\cdot \Omega^{q-1}_X\ra \Omega^q_X)
\setminus \Omega^q_X(-D),
$$
$$f^rh^td(g^s) \wedge \omega\in \Omega^q_{(X,-D)}\setminus 
W_1[{\rm GK}]^q_{(X,-D)},$$
$$f^rg^sh^{t-1}dh \wedge \omega\in 
W_1[{\rm JSZ}]^q_{(X,-D)}\setminus 
\Omega^q_{(X,-D)}.$$
\item 
Both inclusions in \eqref{Comp2} can be strict, this answers \cite[Remark 11.7(1)]{RamFil1}.
If we take the example from (1), then 
$$\dlog(1+f^rg^sh^t)\wedge\dlog h \in
W_1[{\rm JSZ}]^2_{(X,-D),\log}.\setminus\Omega^2_{(X,-D),\log},$$
$$\dlog(1+f^rg^sh^t)\wedge\dlog g \in
\Omega^2_{(X,-D),\log}\setminus W_1[{\rm GK}]^2_{(X,-D),\log}.$$
We can also see the second strict inclusion of \eqref{Comp2} directly from \Cref{thm1}. 
\end{enumerate}
\end{proof}
\begin{remark}
We make a comment on the local freeness of these modules for $n=1$.
While $\Omega^q_{(X,-D)}$ and $W_1[{\rm JSZ}]^q_{(X,-D)}$ are both locally free, $W_1[{\rm GK}]^q_{(X,-D)}$ is not a Cohen-Macaulay $\cO_X$-module on $X_\Zar$.

We can use \Cref{lemAppHartog} to see this. In fact, if it is Cohen-Macaulay, then any rational section which is integral at all codimension 1 points is an integral section.
Take any  $\alpha\in  \Gamma(X, \Omega^q_{(X,-D)})\setminus \Gamma(X, W_1[{\rm GK}]^q_{(X,-D)})$.
In particular, $\alpha\in \Omega^q_{(X,-D),x}$ for all $x\in X^{(1)}$.
Note that the codimension $1$ stalks of $\Omega^q_{(X,-D)}$ and of $W_1[{\rm GK}]^q_{(X,-D)}$ are the same.
We have hence find a section $\alpha$ which lies in $W_1[{\rm GK}]^q_{(X,-D),x}$ for all $x\in X^{(1)}$ but cannot be extended to a global section. 
Therefore $W_1[{\rm GK}]^q_{(X,-D)}$ is not Cohen-Macaulay.
\end{remark}

\section{{Proof of \Cref{MainProp}}}
This section is completely devoted to the proof of \Cref{MainProp}.
The proof is inspired by the proofs of \cite[Proposition 1.2.3]{JSZ} and of \cite[Proposition 1]{Kato82}.
The direction ``$\supset$'' of \Cref{MainPropeq0} is straightforward.
It remains to show the reverse direction.

To this end, we need a small reduction step.
By Artin approximation \cite[Theorem 1.10]{Artin} and the multivariate Hensel's Lemma, any finitely generated $R$-subalgebra $S$ of the completion $\hat R$ of $R$  admits an $R$-homomorphism $S\ra R$. 
Hence for any element from $G^{\ul r}_{R, \log}$, if we can show that it is of the form
\eq{MainPropRedeq}{\sum_j^{\text{finite}}\dlog x_{1,j}\wedge \dlog x_{2,j}\wedge\cdots \wedge\dlog x_{q,j}}
for some
$x_{1,j}\in 
1+(T_1^{\tilde r_1}\cdots  T_g^{\tilde r_g})\cdot \hat R,\,
x_{2,j},\cdots,x_{q,j} \in \hat R[\frac{1}{T_1\cdots T_f}]^\times$,
then it lies in the right hand side of \eqref{MainPropeq0} with $R$ replaced by the finitely generated $R$-subalgebra $S$ of $\hat R$ generated by these finitely many $x_{i,j}$.
With the help of the $R$-homomorphism $S\ra R$, 
we deduce that \eqref{MainPropRedeq} lies in the right hand side of \eqref{MainPropeq0}.
Hence it suffices to show ``$\subset$" with $R$ replaced by $\hat R.$
Moreover, by Cohen structure theorem \cite[\href{https://stacks.math.columbia.edu/tag/032A}{Tag 032A}]{stacks-project}, $\hat R$ is of the form
$$k'[[T_1,\cdots, T_N]]/I$$ 
for some $N\ge d$ and some ideal $I$, 
where $T_1,\cdots,T_d$ in $R$ are the images of the first $d$ variables (which we denote by the same symbols) in $k'[[T_1,\cdots, T_N]]$ via the quotient map.
Hence it suffices to show ``$\subset$" with $R$ replaced by 
$k'[[T_1,\cdots, T_N]]$.

Now we prove the ``$\subset$" direction of \Cref{MainPropeq0} for $R=k'[[T_1,\cdots,T_N]]$.
We first setup some notations.
For each $0\le i\le N$, write 
$$R_i=k'[[T_1,\cdots,\hat T_i,\cdots, T_N]]$$ 
so that $R=R_i[[T_i]]$.
Denote 
$$\tilde A=\Div(T_1\cdots T_f), \quad
\tilde B=\Div(T_1^{\tilde r_1}\cdots  T_g^{\tilde r_g}),$$
and for each $i\in [1,g]$, set
$$\tilde A_i=\Div(T_1\cdots\hat T_i\cdots T_f), \quad
\tilde B_i=\Div(T_1^{\tilde r_1}\cdots \hat T_i^{\tilde r_i}\cdots T_g^{\tilde r_g}),$$
$$ A_i=\Div(T_1\cdots \hat T_i\cdots T_f), \quad
B_i=\Div(T_1^{r_1}\cdots  \hat T_i^{r_i}\cdots T_g^{r_g}).$$
For each $i\in [g+1,N]$, set $A_i=A, B_i=B$ and similarly $\tilde A_i=\tilde A$ and $\tilde B_i=\tilde B$.
(We use $\hat{(-)}$ to mean that the corresponding element does not appear in the string.)
Let $I_i^q$ be the set of strictly increasing functions $\{1,\cdots,q\}\ra \{1,\cdots,\hat i,\cdots , N\}$. 
For $s\in I_i^q$, set 
$$e(s):=\max_{1\le j\le q}\{s(j)\mid s(j)\le e \},
\quad
f(s):=\max_{1\le j\le q}\{s(j)\mid s(j)\le f\}.$$
Denote 
$$\omega_{i,s}= \dlog T_{s(1)}\wedge\cdots \wedge \dlog T_{e(s)}\wedge 
\dlog (1+T_{e(s)+1})\wedge \cdots\wedge \dlog (1+T_{s(q)})
\in \Omega^q_{R_i}(\log A_i),$$
$$\tilde \omega_{i,s}= \dlog T_{s(1)}\wedge\cdots \wedge \dlog T_{f(s)}\wedge 
\dlog (1+T_{f(s)+1})\wedge \cdots\wedge \dlog (1+T_{s(q)})
\in \Omega^q_{R_i}(\log \tilde A_i).$$
Then $\{\omega_{i,s}\mid s\in I_i^q\}$ form a basis of $\Omega^q_{R_i}(\log A_i)$ as a free $R_i$-module.
Similarly, $\{\tilde \omega_{i,s}\mid s\in I_i^q\}$ form a basis of $\Omega^q_{R_i}(\log \tilde A_i)$ as a free $R_i$-module.

Consider the isomorphism of $R_i$-modules (cf. \cite[Claim 1.2.4]{JSZ})
\begin{align}\label{thm1eq1}
\Big(
R\otimes_{R_i}\Omega^q_{R_i}(\log A_i)
\Big)
\oplus
\Big(
R\otimes_{R_i}\Omega^{q-1}_{R_i}(\log A_i)
\Big)
&\xra{\simeq} 
\Omega^q_R(\log A),
\\
(a\otimes w, b\otimes v)
&\mapsto
\begin{cases}
aw+bv\wedge \dlog T_i,
&\text{if $1\le i\le e$;}
\\
aw+bv\wedge \dlog (1+T_i),
&\text{if $e+1\le i\le N$.}
\end{cases}
\nonumber
\end{align}
For any $i\in [1,N]$ and $l\ge 0$, define
$$
V_i^l=\text{ the image via \eqref{thm1eq1} of }
(T^l_i) \otimes_{R_i} \Omega^q_{R_i}(\log A_i)
\oplus
(T^l_i) \otimes_{R_i} \Omega^{q-1}_{R_i}(\log A_i),
$$
We will refer to elements in $V_i^l$ as the \textit{homogeneous degree $l$ part (in the direction $i$)}.
For each $i$, $\{V_i^l\}_{l\ge 0}$ form a decreasing filtration for $\Omega^q_{R}(\log A)$ which is exhaustive. 
In particular, we have
\eq{GV1}{G_R^{\ul r}\subset V_i^{r_i},
}
\eq{GV2}{
G^{(r_1, \dots, r_i+1, \dots, r_N)}_R
=G^{\underline{r}}_R \cap V_i^{r_i+1}}
for any $i$.
Recall that $\Omega^q_R(\log A)_{\log}:=
\Omega^q_R(\log A)\cap \Omega^q_{\Frac R,\log}$.
Set
$$ 
U_i^l:=V_i^l\cap \Omega^q_R(\log A)_{\log}.$$
With these notations, we have the following easy description of the graded piece of $G^{\ul r}_{R}$:
\begin{lem}
\label{lem1}
Fix $i\in [1,N]$. 
Then \eqref{thm1eq1} induces the following isomorphism of $R_i$-modules
\begin{align*}
(T_1^{r_1}\cdots\hat T_i^{r_i}\cdots T_N^{r_N})\cdot
\Omega^q_{R_i}(\log A_i)\,\oplus \,
&
(T_1^{r_1}\cdots\hat T_i^{r_i}\cdots T_N^{r_N})\cdot
\Omega^{q-1}_{R_i}(\log A_i)
\\
&
\xra{\simeq} 
G^{\ul r}_{R}/G^{(r_1,\cdots,r_{i-1},
r_i+1,
r_{i+1},\cdots, r_N)}_{R}
\\ 
(w,v)&\mapsto 
\begin{cases}
T_i^{r_i}w+ T_i^{r_i}v\wedge \dlog T_i
&i\in [1,e],
\\
T_i^{r_i}w+T_i^{r_i}v\wedge \dlog (1+T_i)
&
i\in[e+1,N].
\end{cases}
\end{align*}
\end{lem}

In this section, we will repeatedly use the following fact below: 
\eq{modUvsmodV}{\text{\textit{if $C,D\in \Omega^q_R(\log A)_{\log}$ and $C\equiv D\mod V_i^{h}$, then $C\equiv D \mod U_i^{h}$.}}}

\begin{lem}
\label{lem2}
Let $h>0$ be a positive integer.
\begin{enumerate}
\item 
\label{lem2(1)}
If $i\in [1,e]$, 
any $q$-form in $U_i^{h}/U_i^{h+1}$
is of the form
\begin{align*}
\sum_{s=1}^{\rm finite}
\dlog(1+a_{i,s}T_i^{h})\wedge  w_{i,s}
+
\sum_{t=1}^{\rm finite}
\dlog(1+b_{i,t}T_i^{h})\wedge  w'_{i,t} \wedge \dlog T_i
\end{align*}
where $a_{i,s},b_{i,t}\in R_i, w_{i,s}$, and
$w_{i,s}\in \Omega^{q-1}_{R_i}(\log A_i)_{\log}$, $w'_{i,t}\in \Omega^{q-2}_{R_i}(\log A_i)_{\log}$.
\item 
\label{lem2(2)}
If $i\in [e+1,g]$, $p\mid h$, then any $q$-form in 
$U_i^{h}/U_i^{h+1}$ is of the  form
$$
\sum_{s=1}^{\rm finite} \left(
\dlog(1+a'_{i,s}T_i^{h+1})\wedge w_{i,s}
+
\dlog(1+a_{i,s}T_i^{h})\wedge w_{i,s}
\right)
+
\sum_{t=1}^{\rm finite}
\dlog(1+b_{i,t}T_i^{h})\wedge w'_{i,t}\wedge \dlog (1+T_i),$$
where $a'_{i,s},a_{i,s},b_{i,t}\in R_i$, and
$w_{i,s}\in \Omega^{q-1}_{R_i}(\log A_i)_{\log}$, $w'_{i,t}\in \Omega^{q-2}_{R_i}(\log A_i)_{\log}$.
\item 
\label{lem2(3)}
If $i\in [e+1,g]$, $p\nmid h(h+1)$, then any $q$-form in 
$U_i^{h}/U_i^{h+1}$ is of the  form
$$
\sum_{s=1}^{\rm finite}
\dlog(1+a_{i,s}T_i^{h+1})\wedge w_{i,s}
+
\sum_{t=1}^{\rm finite}
\dlog(1+b_{i,t}T_i^{h})\wedge w'_{i,t}\wedge \dlog (1+T_i),$$
where $a_{i,s},b_{i,t}\in R_i$, and
$w_{i,s}\in \Omega^{q-1}_{R_i}(\log A_i)_{\log}$, $w'_{i,t}\in \Omega^{q-2}_{R_i}(\log A_i)_{\log}$.
\item 
\label{lem2(4)}
If $i\in [e+1,g]$, $p\mid (h+1)$, then any $q$-form in 
$v\in U_i^{h}/U_i^{h+1}$ is of the  form
$$
\sum_{t=1}^{\rm finite}
\dlog(1+b_{i,t}T_i^{h})\wedge w'_{i,t}\wedge \dlog (1+T_i),$$
where $b_{i,t}\in R_i$, and
$w'_{i,t}\in \Omega^{q-2}_{R_i}(\log A_i)_{\log}$.
\end{enumerate}
\end{lem}
\begin{proof}
Every element in $R=R_i[[T_i]]$ which is a unit can be written uniquely as an infinite product
$$u\prod_{l\ge 0}(1-a_lT_i^l),\quad u\in k'^\times, a_l\in R_i.$$
Therefore $\Omega^1_{R,\log}$ is generated by convergent power series of
$$\dlog (1+a_{l}T_i^l)$$
with $a_{l}\in R_i$ and $l\ge 0$.
Moreover, for any $j\in [1,N]$, the natural map
$$R^\times \oplus \Z^j\ra R[\frac{1}{T_1\cdots T_j}]^\times, \quad
(f,(g_1,\cdots ,g_j))\mapsto f\cdot T_1^{g_1}\cdots T_j^{g_j}
$$
is an isomorphism.
In particular, $\Omega^1_{R[\frac{1}{T_1\cdots T_N}],\log}$ is generated by convergent power series of
$$\dlog (1+a_{l}T_i^l), \,
\dlog T_1,\,\cdots,\,\dlog T_N$$
with $a_{l}\in R_i$, $\l\ge 0$.
An element of $\Omega^1_{R[\frac{1}{T_1\cdots T_N}],\log}$ lies in
$\Omega^1_R(\log A)_{\log}$ if  its residue in terms of $T_{i}$  lies in $\Omega^{q-1}_{R_i}(\log A_i)$, for all $i=e+1,\cdots, N$.
Therefore 
$\Omega^1_R(\log A)_{\log}$ 
is generated by convergent power series of
\begin{align}
\label{lem3:eq0}
\dlog(1+a_lT_i^l),
\quad (a_l\in R_i, l\ge 0)
\\
\dlog T_1,\,\cdots,\,\dlog T_e.
\end{align}
(Note that $\dlog (1+T_{e+1}),\cdots, \dlog (1+T_N)$ are among these generators.)
\begin{Observation}
\label{obs1}
Let $w\in \Omega^{q-1}_{R_i}(\log A_i), w'\in \Omega^{q-2}_{R_i}(\log A_i)$, and $a\in R_i$. 
\begin{enumerate}
    \item 

If $i\in [1,e]$ and $h\ge 1$, then 
$$\dlog(1+aT_i^{h})\wedge w =\frac{haT_i^{h}\dlog T_i
+T_i^{h}da
}{1+aT_i^{h}}\wedge w
\in V_i^{h},
$$
$$
\dlog(1+aT_i^{h})\wedge w'\wedge\dlog T_i
=\frac{T_i^{h}da
}{1+aT_i^{h}}\wedge w' \wedge\dlog T_i
\in V_i^{h}.$$
\item 
If $i\in[e+1,g]$, $p\mid h\ge 1$, then
$$\dlog(1+aT_i^{h})\wedge w
=\frac{T_i^{h}da
}{1+aT_i^{h}}\wedge w
\in V_i^{h},$$
$$
\dlog(1+aT_i^{h})\wedge w'\wedge\dlog (1+T_i)
=\frac{T_i^{h}da
}{1+aT_i^{h}}\wedge w' \wedge\dlog (1+T_i)
\in V_i^{h}
.$$
\item 
But if $i\in[e+1,g]$ and $p\nmid h\ge 1$, then
$$\dlog(1+aT_i^{h})\wedge w=
\frac{haT_i^{h-1}(1+T_i)\dlog (1+T_i)+T_i^{h}da
}{1+aT_i^{h}}\wedge w
\in V_i^{h-1},$$
$$
\dlog(1+aT_i^{h})\wedge w'\wedge\dlog (1+T_i)
=\frac{T_i^{h}da
}{1+aT_i^{h}}\wedge w' \wedge\dlog (1+T_i)
\in V_i^{h}
.$$
\end{enumerate}
\end{Observation}

\begin{Observation}\label{obs2}
By \cite[Lemma 2.7(1)]{RS18}, for $l_1,l_2\ge 1$ and $a,b\in R_i$,
$$\{1+aT_i^{l_1},1+bT_i^{l_2}\}=-\{1+\frac{ab}{1+aT_i^{l_1}}T_i^{l_1+l_2},-a(1+bT_i^{l_2})\cdot T_i^{l_1}\}
\quad\text{in $K^M_2(\Frac R)$.}$$
Hence if there are more than one factor of the shape \eqref{lem3:eq0} with $l\ge 1$ in a wedge product, one can always reorganize them so that only one such factor appears.
\end{Observation}

We continue the proof of \Cref{lem2}.
\eqref{lem2(1)} is straightforward.
It remains to prove \eqref{lem2(2)}-\eqref{lem2(4)}.
Let $h_1\in p\Z$ be the unique integer with $h\in[h_1,h_1+p)$.
Write $h=h_0+h_1$ with $h_0\in [0,p)$.
By \Cref{obs1} and \Cref{obs2}, a general element $v\in U^{h_1}_i$, when regarded as an element in the quotient group $U^{h_1}_i/U_i^{h_1+p}$, is of the form
\eq{genformv}{
\sum_{l=0}^{p-1}\sum_{s=1}^{\rm finite}
\dlog(1+a^{(l)}_{i,s}T_i^{h_1+l})\wedge w_{i,s}
+\sum_{l=0}^{p-1}\sum_{t=1}^{\rm finite}
\dlog(1+b^{(l)}_{i,t}T_i^{h_1+l})\wedge w'_{i,t}\wedge \dlog (1+T_i).
}
where all $a^{(l)}_{i,s}, b^{(l)}_{i,t}\in R_i$, and $w_{i,s}\in \Omega^{q-1}_{R_i}(\log A_i)_{\log}$, $w'_{i,t}\in \Omega^{q-2}_{R_i}(\log A_i)_{\log}$.
\eqref{lem2(2)} then follows directly from this expression.
It remains to discuss \eqref{lem2(3)} and \eqref{lem2(4)}.
It will be sufficient to prove the following claim.

\begin{claim}
\label{lem2(34)claim}
If $i\in [e+1,g]$, $p\nmid h$, then any $q$-form in 
$U_i^{h}/U_i^{h+1}$ is of the  form
$$
\sum_{s=1}^{\rm finite}
\dlog(1+a_{i,s}T_i^{h+1})\wedge w_{i,s}
+
\sum_{t=1}^{\rm finite}
\dlog(1+b_{i,t}T_i^{h})\wedge w'_{i,t}\wedge \dlog (1+T_i),$$
where $a_{i,s},b_{i,t}\in R_i$, and
$w_{i,s}\in \Omega^{q-1}_{R_i}(\log A_i)_{\log}$, $w'_{i,t}\in \Omega^{q-2}_{R_i}(\log A_i)_{\log}$.
\end{claim}

\textbf{Case 1. $h_1\ge p$.}
In this case, we have
$$(a^{(l)}_{i,s}T_i^{h_1+l})^{l'}d(a^{(l)}_{i,s}T_i^{h_1+l})\wedge  w_{i,s}
\in 
V_i^{h_1+p},$$
$$(b^{(l)}_{i,t}T_i^{h_1+l})^{l'}\cdot T_i^{h_1+l}db^{(l)}_{i,t}\wedge w_{i,t}\wedge \dlog (1+T_i)
\in V_i^{h_1+p}$$
for all $l\in [0,p-1]$, $l'\ge 1$.
These two relations imply that
$$\dlog(1+a^{(l)}_{i,s}T_i^{h_1+l})\wedge w_{i,s}
\equiv 
d(a^{(l)}_{i,s}T_i^{h_1+l})\wedge w_{i,s} 
\mod V_i^{h_1+p},$$
$$\dlog(1+b^{(l)}_{i,t}T_i^{h_1+l})\wedge w_{i,t}\wedge \dlog (1+T_i)
\equiv 
T_i^{h_1+l} db^{(l)}_{i,t}\wedge w_{i,t}\wedge \dlog (1+T_i)
\mod V_i^{h_1+p}.$$
We will repeatedly use these congrence relations below.
Now let $v\in U^{h_1+h_0}_i/U_i^{h_1+h_0+1}$.
We prove \Cref{lem2(34)claim} by discussing two cases of $h_0$.

\textbf{Subcase 1.1.} $h_0=p-1$.
Let
$v \in U^{h_1+p-1}_i/U_i^{h_1+p}$.
A direct computation from \eqref{genformv} shows
\begin{align*}
v&\equiv 
\sum_{l=0}^{p-1}\sum_{s=1}^{\rm finite}
T_i^{h_1+l}da_{i,s}^{(l)}\wedge w_{i,s}
\\&+
\sum_{l=0}^{p-2}
T_i^{h_1+l}
\left(
\sum_{s=1}^{\rm finite}
(-1)^{q-1}
(h_1+l+1) a_{i,s}^{(l+1)} w_{i,s}
+
\sum_{t=1}^{\rm finite}
\left(\sum_{l'=0}^{l} (-1)^{l-l'}db_{i,t}^{(l')}\right)
\wedge  w_{i,t}
\right)
\wedge dT_i
\\
&+
T_i^{h_1+p-1}
\left(\sum_{t=1}^{\rm finite}(db_{i,t}^{(p-1)}-db_{i,t}^{(p-2)}+\cdots +(-1)^{p-1}db_{i,t}^{(0)})\wedge  w_{i,t}\right)
\wedge dT_i
\qquad \mod V_i^{h_1+p}.
\end{align*}
The condition $v\in V_i^{h_1+p-1}$ forces
\eq{lem2(34)claim1.1-1}{\sum_{s=1}^{\rm finite}  da_{i,s}^{(l)}\wedge w_{i,s}=0,
\quad
\text{for all $l\in [0,p-2]$; and}}
\eq{lem2(34)claim1.1-2}{\sum_{s=1}^{\rm finite}
(-1)^{q-1} (h_1+l+1)a_{i,s}^{(l+1)} w_{i,s}
+
\sum_{t=1}^{\rm finite}
\left(\sum_{l'=0}^{l} (-1)^{l-l'}db_{i,t}^{(l')}\right)
\wedge  w_{i,t}
=0,
\quad
\text{for all $l\in [0,p-2]$.}}
Taking $l=p-2$ in the second vanishing condition above and differentiating it, we get
$$\sum_{s=1}^{\rm finite} 
da_{i,s}^{(p-1)}\wedge w_{i,s}
=0.$$
That is, \eqref{lem2(34)claim1.1-1} holds for all $l\in [0,p-1]$.
Therefore
\begin{align*}
v&\equiv 
T_i^{h_1+p-1}
\left(
\sum_{t=1}^{\rm finite}
(db_{i,t}^{(p-1)}-db_{i,t}^{(p-2)}+\cdots +(-1)^{p-1}db_{i,t}^{(0)})\wedge  w_{i,t}\right)
\wedge dT_i
\mod V_i^{h_1+p}
\\
&\equiv 
\sum_{t=1}^{\rm finite}
\dlog(1+
\left(\sum_{l=0}^{p-1} (-1)^{p-1-l}b_{i,t}^{(l)}\right)
T_i^{h_1+p-1})\wedge w_{i,t}\wedge \dlog (1+T_i)
\mod V_i^{h_1+p}.
\end{align*}
Since both sides of the congruence equation above lie in $U_i^{h_1+p-1}$, their difference lies in $U_i^{h_1+p-1}\cap V_i^{h_1+p}=U_i^{h_1+p}$ by \eqref{modUvsmodV}.
In other words, any general element $v\in U^{h_1+p-1}_i/U_i^{h_1+p}$ is of the form
$$
\sum_{t=1}^{\rm finite}
\dlog(1+b_{i,t}T_i^{h_1+p-1})\wedge w_{i,t}\wedge \dlog (1+T_i)
$$
for some $b_{i,t}\in R_i$.
In particular, \Cref{lem2(34)claim} holds in this subcase.

\textbf{Subcase 1.2.} $h_0\in [0,p-1).$
Let $v\in U^{h_1+h_0}_i/U_i^{h_1+h_0+1}$.
Following the general form \eqref{genformv},  we have
\begin{align*}
v
&\equiv 
\sum_{l=0}^{h_0+1}\sum_{s=1}^{\rm finite}
\dlog(1+a^{(l)}_{i,s}T_i^{h_1+l})\wedge w_{i,s}
+\sum_{l=0}^{h_0}\sum_{t=1}^{\rm finite}
\dlog(1+b^{(l)}_{i,t}T_i^{h_1+l})\wedge w_{i,t}\wedge \dlog (1+T_i)
\\
&\mod V_i^{h_1+h_0+1}.
\end{align*}
A direct computation shows
\begin{align*}
v&\equiv 
\sum_{l=0}^{h_0}\sum_{s=1}^{\rm finite}
T_i^{h_1+l}da_{i,s}^{(l)}\wedge w_{i,s}
\\&+
\sum_{l=0}^{h_0-1}
T_i^{h_1+l}
\left(
\sum_{s=1}^{\rm finite} 
(-1)^{q-1} (h_1+l+1)
a_{i,s}^{(l+1)} w_{i,s}
+
\sum_{t=1}^{\rm finite}
\left(\sum_{l'=0}^{l} (-1)^{l-l'} db_{i,t}^{(l')}\right)
\wedge  w_{i,t}
\right)
\wedge dT_i
\\
&+
T_i^{h_1+h_0}
\left(
\sum_{s=1}^{\rm finite} 
(-1)^{q-1} (h_1+h_0+1)
a_{i,s}^{(h_0)} w_{i,s}
+
\sum_{t=1}^{\rm finite}
\left(\sum_{l'=0}^{h_0} (-1)^{h_0-l'} db_{i,t}^{(l')}\right)
\wedge  w_{i,t}
\right)
\wedge dT_i
\\
&
\mod V_i^{h_1+h_0+1}.
\end{align*}
The condition $v\in V_i^{h_1+h_0}$ forces
\eq{lem2(34)claim1.2-1}{\sum_{s=1}^{\rm finite}da_{i,s}^{(l)}\wedge w_{i,s}=0,
\quad
\text{for all $l\in [0,h_0-1]$; and}}
\eq{lem2(34)claim1.2-2}{\sum_{s=1}^{\rm finite} (-1)^{q-1} (h_1+l+1) a_{i,s}^{(l+1)} w_{i,s}
+
\sum_{t=1}^{\rm finite}
\left(\sum_{l'=0}^{l} (-1)^{l-l'} db_{i,t}^{(l')}\right)
\wedge  w_{i,t}
=0,
\quad
\text{for all $l\in [0,h_0-1]$.}}
Taking $l=h_0-1$ in the second vanishing condition above and differentiating it, we get
$$\sum_{s=1}^{\rm finite}
da_{i,s}^{(h_0-1)}\wedge w_{i,s}=0.$$
That is, \eqref{lem2(34)claim1.2-1} holds for all $l\in [0,h_0]$.
Therefore
\begin{align*}
v&\equiv 
T_i^{h_1+h_0}
\left(
\sum_{s=1}^{\rm finite} 
(-1)^{q-1} (h_1+h_0+1) a_{i,s}^{(h)} w_{i,s}
+
\sum_{t=1}^{\rm finite}
\left(\sum_{l=0}^{h_0} (-1)^{h_0-l} db_{i,t}^{(l)}\right)
\wedge  w_{i,t}
\right)
\wedge dT_i
\\
&\quad
\mod V_i^{h_1+h_0+1}
\\
&\equiv 
\sum_{s=1}^{\rm finite}
\dlog(1+a^{(h_0)}_{i,s}T_i^{h_1+h_0+1})\wedge w_{i,s}
\\
&+\sum_{t=1}^{\rm finite}
\dlog(1+
\left(\sum_{l=0}^{h_0} (-1)^{h_0-l}b_{i,t}^{(l)}\right)
T_i^{h_1+h_0})\wedge w_{i,t}\wedge \dlog (1+T_i)
\mod V_i^{h_1+h_0+1}.
\end{align*}
Since both sides of the congruence equation above lies in $U_i^{h_1+h_0}$, their subtraction lies in $U_i^{h_1+h_0}\cap V_i^{h_1+h_0+1}=U_i^{h_1+h_0+1}$ by \eqref{modUvsmodV}.
This finishes the proof of \Cref{lem2(34)claim} in subcase 1.2.

\textbf{Case 2. $h\in [1,p)$.}
A general element $U^{0}_i/U_i^{p}$ is of the form
\eq{genformv(1)}{
\sum_{l=0}^{p-1}\sum_{s=1}^{\rm finite}
\dlog(1+a^{(l)}_{i,s}T_i^{l})\wedge w_{i,s}
+\sum_{l=0}^{p-1}\sum_{t=1}^{\rm finite}
\dlog(1+b^{(l)}_{i,t}T_i^{l})\wedge w_{i,t}\wedge \dlog (1+T_i).
}
(Due to formatting reasons, below we will freely omit some indices which should be clear from the context.)
We prove the general form of an element $v\in U^{h}_i/U_i^{h+1}$.

\textbf{Subcase 2.1.} $h=p-1$.
Let $v \in U^{p-1}_i/U^p_i$. Equation \eqref{genformv(1)} implies
\begin{align*}
v&=
\sum_s \dlog(1+a^{(0)})\wedge  w_{i,s}
\\&+
\sum_{l=1}^{p-1}\sum_s 
\left(
\sum_{l'\ge 0}^{\infty}
(-a^{(l)}T^l)^{l'}\cdot (T^lda^{(l)}+la^{(l)}T^{l-1}dT)
\right)\wedge  w_{i,s}
\\&+
\sum_t \dlog(1+b^{(0)})\wedge  w_{i,t} \wedge
\left(\sum_{l''\ge 0}^{\infty} (-T)^{l''}\right)dT
\\&+
\sum_{l=1}^{p-1}\sum_t 
\left( 
\sum_{l'\ge 0}^{\infty}
(-b^{(l)}T^l)^{l'}\cdot T^ldb^{(l)}
\right)
\wedge  w_{i,t}
\wedge 
\left(\sum_{l''\ge 0}^{\infty} (-T)^{l''}\right)dT
\\
&\equiv
\sum_s \dlog(1+a^{(0)})\wedge  w_{i,s}
\\&+
\sum_{l=1}^{p-1} \sum_s 
T^{l}\left(
\sum_{l'\mid l} (-a^{(l')})^{l/l'-1}da^{(l')}
\right)\wedge w_{i,s}
\\&+
\sum_{l=0}^{p-1}  
T^l\Bigg(
\sum_s
\Bigg(
\sum_{l'\mid (l+1)} l' a^{(l')}\cdot  (-a^{(l')})^{(l+1)/l'-1}
\Bigg)\cdot  w_{i,s}
\\
&\quad
+
\sum_t
\Bigg(
(-1)^l\dlog(1+b^{(0)})
+
\sum_{l'=1}^l \sum_{
\stackrel{l''\ge 0}{ l'(l''+1)\le l}}
(-1)^{l-l'(l''+1)}(-b^{(l')})^{l''}db^{(l')}
\Bigg)\wedge w_{i,t}
\Bigg)
\wedge dT
\\&
\mod V_i^p.
\end{align*}
The condition that $v\in V_i^{p-1}$ forces
\eq{genformv(1)1}{
\sum_s \dlog(1+a^{(0)})\wedge  w_{i,s}=0;
}
\eq{genformv(1)2}{
\sum_s 
\left(
\sum_{l'\mid l} (-a^{(l')})^{l/l'-1}da^{(l')}
\right)\wedge w_{i,s}=0
\quad
\text{for all $l\in[1,p-2]$;}
}
\begin{align}
\label{genformv(1)3}
&\sum_s 
\Bigg(
\sum_{l'\mid (l+1)}
(-1)^{(l+1)/l'-1} l'  (a^{(l')})^{(l+1)/l'}
\Bigg)\cdot  w_{i,s}
\\
+&
\sum_t 
\Bigg(
(-1)^l\dlog(1+b^{(0)})
+
\sum_{l'=1}^l \sum_{
\stackrel{l''\ge 0}{ l'(l''+1)\le l}}
(-1)^{l-l'(l''+1)}(-b^{(l')})^{l''}db^{(l')}
\Bigg)\wedge w_{i,t}
=0
\nonumber
\\
&\text{for all $l\in[0,p-2]$.}
\nonumber
\end{align}
Taking $l=p-2$ in
\eqref{genformv(1)3}, differentiating it and dividing it by  $(p-1)$, we deduce that
$$\sum_s
\left(
\sum_{l'\mid (p-1)} (-a^{(l')})^{(p-1)/l'-1}da^{(l')}
\right)\wedge w_{i,s}=0.
$$
That is, \eqref{genformv(1)2} holds for all $l\in [1,p-1]$.
As a result,
\begin{align}
\label{genformv(1)4}
v&\equiv
T^{p-1}\Bigg(
\sum_s
\Bigg(
\sum_{l'\mid p} l' a^{(l')}\cdot  (-a^{(l')})^{p/l'-1}
\Bigg)\cdot w_{i,s}
+
\sum_t
\Bigg(
(-1)^{p-1}\dlog(1+b^{(0)})
\\
&\quad 
+
\sum_{l'=1}^{p-1} \sum_{
\stackrel{l''\ge 0}{ l'(l''+1)\le p-1}}
(-1)^{p-1-l'(l''+1)+l''}(b^{(l')})^{l''}db^{(l')}
\Bigg)\wedge w_{i,t}
\Bigg)
\wedge dT
\mod V_i^p
\nonumber
\\
&\equiv
T^{p-1}
\sum_s
\Bigg(
(-1)^{p-1}\cdot  (a^{(1)})^{p}
\Bigg)
\cdot  w_{i,s}
\wedge dT
\nonumber
\\
&\quad 
+
T^{p-1}\Bigg(
\sum_t
\Bigg(
(-1)^{p-1}\dlog(1+b^{(0)})
\nonumber
\\&
\quad 
+
\sum_{l'=1}^{p-1} \sum_{
\stackrel{l''\ge 0}{ l'(l''+1)\le p-1}}
(-1)^{p-1-l'(l''+1)+l''}(b^{(l')})^{l''}db^{(l')}
\Bigg)\wedge w_{i,t}
\Bigg)
\wedge \dlog(1+T)
\nonumber
\mod V_i^{p}
\\
&\equiv
\sum_s
(-1)^{p-1} 
C^{-1}( w_{i,s}
\wedge a^{(1)}dT)
\nonumber
\\
&\quad 
+
\sum_t (-1)^{p-1}
C^{-1}(\dlog(1+b^{(0)})\wedge w_{i,t}\wedge dT)
\nonumber
\\&
+
\sum_t 
\sum_{l'\ge 1}^{p-1} 
\sum_{l''\ge 0}^{\lfloor \frac{p-1}{l'}\rfloor -1} (-1)^{p-1}
\frac{(-1)^{l'l''+l'+l''}}{l''+1}
\dlog(1+(b^{(l')})^{l''+1}T^{p-1})\wedge w_{i,t}\wedge \dlog(1+T)
\mod V_i^{p}.
\nonumber
\end{align} 
In the second congruence relation, we have used that $p a^{(p)}=0$ in $R_i$.
Let 
\begin{center} 
$\sB_n:=C^{-(n-1)}d\Omega_X^{q-1}(\log A)$ \quad 
for any $n\ge 1$.
\end{center}
Note that
\begin{align} 
C^{-1}(a^{(1)}dT\wedge  w_{i,s})
&\equiv
-C^{-1}(Td(a^{(1)} w_{i,s}))
\mod \sB_2
\label{modB2}
\\
&\equiv 0 \mod V_i^p.
\nonumber
\end{align}
In the meanwhile,
\eq{modB21}{
C^{-1}(\dlog(1+b^{(0)})\wedge w_{i,t}\wedge dT)
=(-1)^{q-1}C^{-1}(d(T\dlog(1+b^{(0)})\wedge w_{i,t}))
\in \sB_2.}
We write $$\sB_\infty:=\bigcup_{n\ge 1} \sB_n.$$
By \Cref{sesF-1Binf}, 
the natural map
\eq{logmodBinf}{\Omega^q_R(\log A)_{\log} \ra \Omega^q_R(\log A)/\sB_{\infty}}
is an injection.
In particular, 
$$\Omega^q_{R}(\log A)_{\log}\simeq 
{\rm Image}\left(\Omega^q_{R}(\log A)_{\log}\hra \frac{\Omega^q_{R}(\log A)}{\sB_\infty}\right)=
{\rm Image}\left(\Omega^q_{R}(\log A)_{\log}+\sB_\infty\ra \frac{\Omega^q_{R}(\log A)}{\sB_\infty}\right).$$
Moreover, for any $l$, since $\Omega^q_{R}(\log A)_{\log}\cap (V_i^l+\sB_\infty)$ is the preimage of the image of $\Omega^q_{R}(\log A)_{\log}\cap V_i^l$ under the injective map \eqref{logmodBinf}, we have
\eq{modVlismodVlplusB}{\Omega^q_{R}(\log A)_{\log}\cap V_i^l=\Omega^q_{R}(\log A)_{\log}\cap (V_i^l+\sB_\infty).}
Now \eqref{genformv(1)4}, \eqref{modB2}, \eqref{modB21} and \eqref{modVlismodVlplusB} together imply
$$v\equiv 
\sum_t \sum_{l'\ge 1}^{p-1} \sum_{l''\ge 0}^{\frac{p-1}{l'}-1} (-1)^{p-1}
\frac{(-1)^{l'l''+l'+l''}}{l''+1}
\dlog(1+(b^{(l')})^{l''+1}T_i^{p-1})\wedge w_{i,t}\wedge \dlog(1+T_i)
\mod U_i^p.$$
Taking into account that
\begin{align*} 
\dlog(1+bT_i^{p-1})&\wedge w_{i,t}\wedge \dlog(1+T_i)+
 \dlog(1+b'T_i^{p-1})\wedge w_{i,t}\wedge \dlog(1+T_i)
\\
&\equiv
 \dlog(1+(b+b')T_i^{p-1})\wedge w_{i,t}\wedge \dlog(1+T_i)
\mod U_i^p
 \end{align*}
for any $b,b'\in R_i$, we deduce that any $v\in  U^{p-1}_i/U_i^{p}$ is of the form
$$\sum_t \dlog(1+bT_i^{p-1})\wedge w_{i,t}\wedge \dlog(1+T_i),\quad b\in R_i.$$
In particular, \Cref{lem2(34)claim} holds in this case.

\textbf{Subcase 2.2.} $h \in[1, p-1)$. 
Let $v \in U^h_i / U^{h+1}_i$. By \eqref{genformv(1)}, since $v \in U^h_i$, we can truncate the expansion at the relevant degrees to write
\[
v \equiv \sum_{l=0}^{h+1} \sum_s \mathrm{dlog}(1 + a_{i,s}^{(l)}T_i^l) \wedge w_{i,s} + \sum_{l=0}^h \sum_t \mathrm{dlog}(1 + b_{i,t}^{(l)}T_i^l) \wedge w_{i,t} \wedge \mathrm{dlog}(1 + T_i) \pmod{U^{h+1}_i}.
\]
Similar to the calculations before, 
the condition that $v\in V_i^{p-1}$ forces
\eq{genformv(1)11}{\sum_s\dlog(1+a^{(0)})\wedge  w_{i,s}=0}
\eq{genformv(1)22}{
\sum_s
\left(
\sum_{l'\mid l} (-a^{(l')})^{l/l'-1}da^{(l')}
\right)\wedge w_{i,s}=0
\quad
\text{for all $l\in[1,h-2]$,}
}
\begin{align}
\label{genformv(1)33}
&\sum_s 
\Bigg(
\sum_{l'\mid (l+1)}
(-1)^{(l+1)/l'-1} l'  (a^{(l')})^{(l+1)/l'}
\Bigg)\cdot  w_{i,s}
\\
+&
\sum_t 
\Bigg(
(-1)^l\dlog(1+b^{(0)})
+
\sum_{l'=1}^l \sum_{
\stackrel{l''\ge 0}{ l'(l''+1)\le l}}
(-1)^{l-l'(l''+1)}(-b^{(l')})^{l''}db^{(l')}
\Bigg)\wedge w_{i,t}
=0
\nonumber
\\
&\text{for all $l\in[0,h-2]$.}
\nonumber
\end{align}
Taking $l=h-2$ in
\eqref{genformv(1)33}, differentiating it and dividing the result by  $(h-1)$ 
(note that $1\le h-1\le p-2$), we deduce that
$$\sum_s 
\left(
\sum_{l'\mid (h-1)} (-a^{(l')})^{(p-1)/l'-1}da^{(l')}
\right)\wedge w_{i,s}=0.
$$
That is, \eqref{genformv(1)22} holds for all $l\in [1,h-1]$.
Hence
\begin{align*} 
v&\equiv
T^{h-1}\Bigg(
\sum_s
\Bigg(
\sum_{l'\mid h} l' a^{(l')}\cdot  (-a^{(l')})^{h/l'-1}
\Bigg)\cdot w_{i,s}
+
\sum_t
\Bigg(
(-1)^{h-1}\dlog(1+b^{(0)})
\\
&\quad 
+
\sum_{l'=1}^{h-1} \sum_{
\stackrel{l''\ge 0}{ l'(l''+1)\le h-1}}
(-1)^{h-1-l'(l''+1)+l''}(b^{(l')})^{l''}db^{(l')}
\Bigg)\wedge w_{i,t}
\Bigg)
\wedge dT
\mod V_i^h
\\
&
\equiv
T^{h-1}\Bigg(
\sum_s
\sum_{l'\mid h}
(-1)^{h/l'-1} l'\cdot  (a^{(1)})^{h/l'}
\cdot  w_{i,s}
\Bigg)
\wedge \dlog(1+ T)
\\&\quad +
T^{h-1}\Bigg(
\sum_t
\Bigg(
(-1)^{h-1}\dlog(1+b^{(0)})
\nonumber
\\&\quad +
\sum_{l'=1}^{h-1} \sum_{
\stackrel{l''\ge 0}{ l'(l''+1)\le h-1}}
(-1)^{h-1-l'(l''+1)+l''}(b^{(l')})^{l''}db^{(l')}
\Bigg)\wedge w_{i,t}
\Bigg)
\wedge \dlog(1+T)
\mod V_i^h
\\
&
\equiv
 \sum_s
\sum_{l'\mid h}
\frac{(-1)^{h/l'-1} l'}{h}\cdot  w_{i,s}
\wedge \dlog(1+ (a^{(1)})^{h/l'}T^h)
\\&\quad +
  \sum_t
(-1)^{h-1}\dlog(1+b^{(0)})
\wedge w_{i,t}
\wedge T^{h-1}dT
\nonumber
\\& \quad +
 \sum_t\sum_{l'=1}^{h-1} \sum_{
\stackrel{l''\ge 0}{ l'(l''+1)\le h-1}}
\frac{(-1)^{p-1+l'(l''+1)+l''}}{l''+1}
\dlog(1+(b^{(l')})^{l''+1}T^{h-1})
\wedge w_{i,t}
\wedge \dlog(1+T)
\\
&
\mod V_i^h.
\end{align*} 
Note that
\eq{}{
\dlog(1+b^{(0)})
\wedge w_{i,t}
\wedge T^{h-1}dT
=\frac{(-1)^{q-1}}{h}
d(T^h\dlog(1+b^{(0)})
\wedge w_{i,t})
\in \sB_1.
}
\eqref{modVlismodVlplusB} implies
\begin{align*} 
v&\equiv
 \sum_s
\sum_{l'\mid h}
\frac{(-1)^{h/l'-1} l'}{h}\cdot  w_{i,s}
\wedge \dlog(1+ (a^{(1)})^{h/l'}T^h)
\\
\\&+
 \sum_t\sum_{l'=1}^{h-1} \sum_{
\stackrel{l''\ge 0}{ l'(l''+1)\le h-1}}
\frac{(-1)^{p-1+l'(l''+1)+l''}}{l''+1}
\dlog(1+(b^{(l')})^{l''+1}T^{h-1})
\wedge w_{i,t}
\wedge \dlog(1+T)
\mod U_i^h.
\end{align*}
This last expression  and the additivity of such forms at  the positions
$(a^{(1)})^{h/l'}$ and $(b^{(l')})^{l''+1}$ conclude the proof of Subcase 2.2.
\end{proof}

\begin{lem}
\label{lem3}
Set 
$$\cZ^{q}:=
\Ker(d: \Omega^{q}_{R_i}(\log A_i)\ra \Omega^{q+1}_{R_i}(\log A_i)).
$$
As a convention, we set $\Omega^{-1}=0$, and
$$\Omega^0_{R_i}(\log A_i)=R_i.$$
Let $h\ge 0$ be an integer.
\begin{enumerate}
\item\label{lem3(1} 
If $e>0$, $i\in [1,e]$ and $p\nmid h$, then we have the following isomorphism of abelian groups
\begin{align*}
\rho_{i,h}:
\Omega^{q-1}_{R_i}(\log A_i)
&\xra{\simeq}
U_i^{h}/U_i^{h+1}
\\
\sum_{s\in I_i^{q-1}} a_{i,s}\omega_{i,s}
&\mapsto
    \sum_{s\in I_i^{q-1}} \dlog(1+a_{i,s}T_i^{h})\wedge\omega_{i,s},
\end{align*}
where $a_{i,s}\in R_i$.

Moreover, 
$$0\ra \frac{\Omega^{q-2}_{R_i}(\log A_i)}{\cZ^{q-2}}
\xra{\beta} U_i^{h}/U_i^{h+1}
\xra{\alpha} \frac{\Omega^{q-1}_{R_i}(\log A_i)}{\cZ^{q-1}}
\ra 0 $$
is a short exact sequence, where $\beta$ is defined by
\eq{betaeq}{
\beta(\sum_t b_{i,t}\omega_{i,t})=\sum_{t} \dlog(1+b_{i,t}T_i^{h})\wedge\omega_{i,t}\wedge\dlog T_i
}
and $\alpha$ is defined by the inverse of $\rho_{i,h}$ followed by the projection.
\item \label{lem3(2}
If $e>0$, $i\in[1,e]$ and $p\mid h\ge 1$, then we have the following isomorphism of abelian groups
\begin{align*}
\rho_{i,h}:
\frac{\Omega^{q-1}_{R_i}(\log A_i)}
{\cZ^{q-1}}
\oplus
\frac{\Omega^{q-2}_{R_i}(\log A_i)}
{\cZ^{q-2}}
&\xra{\simeq}
U_i^{h}/U_i^{h+1}
\\
\left(
\sum_{s\in I_i^{q-1}} a_{i,s}\omega_{i,s},
\sum_{t\in I_i^{q-2}} b_{i,t}\omega_{i,t}
\right) 
&\mapsto
\sum_{s\in I_i^{q-1}} 
\dlog(1+a_{i,s}T_i^{h})\wedge\omega_{i,s}
\\
&\quad +
\sum_{t\in I_i^{q-2}} 
\dlog(1+b_{i,t}T_i^{h})\wedge\omega_{i,t}
\wedge \dlog T_i,
\end{align*}
where $a_{i,s}, b_{i,t}\in R_i$.

\item \label{lem3(3}
If $e>0$, $i\in [1,e]$ and $h=0$, then we have the following isomorphism of abelian groups
\begin{align*}
\rho_{i,h}:
\Omega^{q}_{R_i}(\log A_i)_{\log}
\oplus
\Omega^{q-1}_{R_i}(\log A_i)_{\log}
&\xra{\simeq}
U_i^{0}/U_i^{1}
\\
(w,w')&\mapsto w+w'\wedge \dlog T_i.
\end{align*}

\item \label{lem3(4}
If  $i\in [e+1,g]$ and $p\mid h\ge 1$, then we have the following isomorphism of abelian groups
\begin{align*}
\rho_{i,h}:
\Omega^{q-1}_{R_i}(\log A_i)
\oplus
\frac{\Omega^{q-1}_{R_i}(\log A_i)}{\cZ^{q-1}}
&\xra{\simeq}
U_i^{h}/U_i^{h+1},
\\
\Big(
\sum_{s\in I_i^{q-1}} a'_{i,s}\omega_{i,s},
\sum_{s\in I_i^{q-1}} a_{i,s}\omega_{i,s},
\Big)
&\mapsto
\sum_{s\in I_i^{q-1}} \dlog(1+a'_{i,s}T_i^{h+1})\wedge\omega_{i,s}
\nonumber
\\
&\quad +
\sum_{s\in I_i^{q-1}} \dlog(1+a_{i,s}T_i^{h})\wedge\omega_{i,s},
\end{align*}
where $a'_{i,s}, a_{i,s}\in R_i$.

Moreover, 
$$0\ra \frac{\Omega^{q-2}_{R_i}(\log A_i)}{\cZ^{q-2}}
\xra{\beta} U_i^{h}/U_i^{h+1}
\xra{\alpha}
\frac{\Omega^{q-1}_{R_i}(\log A_i)}{\cZ^{q-1}}\oplus \frac{\Omega^{q-1}_{R_i}(\log A_i)}{\cZ^{q-1}}
\ra 0 $$
is a short exact sequence, where $\beta$ is defined by 
\eq{betaeq1}{
\beta(\sum_t b_{i,t}\omega_{i,t})=\sum_{t} \dlog(1+b_{i,t}T_i^{h})\wedge\omega_{i,t}\wedge\dlog (1+T_i),
}
and $\alpha$ is defined by the inverse of $\rho_{i,h}$ followed by the projection.

\item\label{lem3(5}
If $i\in [e+1,g]$ and $p\nmid h(h+1)$,
then we have the following isomorphism of abelian groups
\begin{align*}
\rho_{i,h}:
\Omega^{q-1}_{R_i}(\log A_i)
&\ra
U_i^{h}/U_i^{h+1}
\\
\sum_{s\in I_i^{q-1}} a_{i,s}\omega_{i,s},
&\mapsto
\sum_{s\in I_i^{q-1}} \dlog(1+a_{i,s}T_i^{h+1})\wedge\omega_{i,s}.
\end{align*}
where $a_{i,s}\in R_i$. 

Moreover, 
$$0\ra 
\frac{\Omega^{q-2}_{R_i}(\log A_i)}{\cZ^{q-2}}
\xra{\beta} U_i^{h}/U_i^{h+1}
\xra{\alpha} \frac{\Omega^{q-1}_{R_i}(\log A_i)}{\cZ^{q-1}}
\ra 0 $$
is a short exact sequence, where $\beta$ is defined by \eqref{betaeq1},
and $\alpha$ is defined by the inverse of $\rho_{i,h}$ followed by the projection.

\item\label{lem3(6}
If $i\in [e+1,g]$ and $p\mid (h+1)$, then we have the following isomorphism of abelian groups
\begin{align*}
\rho_{i,h}:
\frac{\Omega^{q-2}_{R_i}(\log A_i)}
{\cZ^{q-2}}
&\ra
U_i^{h}/U_i^{h+1}
\\
\sum_{t\in I_i^{q-2}} b_{i,t}\omega_{i,t}
&\mapsto
\sum_{t\in I_i^{q-2}} 
\dlog(1+b_{i,t}T_i^{h})\wedge\omega_{i,t}
\wedge \dlog (1+T_i).
\end{align*}
where $b_{i,t}\in R_i$.
(In particular, when $q=1$, we have $U_i^{h}=U_i^{h+1}$.)

\item \label{lem3(7}
If $i\in [g+1,N]$ and $h=0$, then we have the following isomorphism of abelian groups
\begin{align*}
\rho_{i,h}:
\Omega^q_{R_i}(\log A_i)_{\log}
\oplus \Omega^{q-1}_{R_i}(\log A_i)_{\log}
&\xra{\simeq}
U_i^{0}/U_i^{1}.
\\
(w,w')
&\mapsto w
+ w'\wedge\dlog(1+T_i).
\end{align*}
\end{enumerate}
\end{lem}
\begin{proof}[Proof of \Cref{lem3}]
The first statements in (1)(2) are already contained in the proof of \cite[p. 224, Lemma 4]{Kato82}. We include the proof here for the convenience of the reader.
\begin{enumerate}[label= (\alph*)]
\item 
If $i\in [1,e]$ and $h\ge 1$, 
consider the map
\begin{align} 
\label{prop:cl:eq1}
\rho_{i,h}:
\Omega^{q-1}_{R_i}(\log A_i)\oplus \Omega^{q-2}_{R_i}(\log A_i)
&\ra
U_i^{h}/U_i^{h+1},
\\
\Big(
\sum_{s\in I_i^{q-1}}a_{i,s}\omega_{i,s},
\sum_{t\in I_i^{q-2}}b_{i,t}\omega_{i,t}
\Big)
&\mapsto
\sum_{s\in I_i^{q-1}} \dlog(1+a_{i,s}T_i^{h})\wedge\omega_{i,s}
\nonumber
\\
&\quad 
+\sum_{t\in I_i^{q-2}}
\dlog(1+b_{i,t}T_i^{h})\wedge\omega_{i,t}\wedge \dlog T_i.
\nonumber
\end{align}
This map is clearly well-defined and is a map of abelian groups.
It is also surjective by \Cref{lem2}\eqref{lem2(1)}.
In fact, take any general element 
\eq{prop:cl:eq12}{\sum_{s\in I_i^{q-1}} 
\dlog(1+a'_{i,s}T_i^{h})\wedge w_{i,s}
+\sum_{t\in I_i^{q-2}}
\dlog(1+b'_{i,t}T_i^{h})\wedge w'_{i,t}\wedge \dlog T_i
}
from $U_i^h$, where
$a'_{i,s},b'_{i,t}\in R_i$, 
$w_{i,s}\in \Omega^{q-1}_{R_i}(\log A_i)_{\log}$, $w'_{i,t}\in \Omega^{q-2}_{R_i}(\log A_i)_{\log}$.
Because $\{\omega_{i,s}\}_s$ form a basis, there exist
$a_{i,s}, b_{i,t}\in R_i$ such that
$$\sum_{s\in I_i^{q-1}} a_{i,s}\omega_{i,s}
=\sum_{s\in I_i^{q-1}} a'_{i,s} w_{i,s},
\quad 
\sum_{t\in I_i^{q-2}} b_{i,t}\omega_{i,t}
=\sum_{t\in I_i^{q-2}} b'_{i,t}  w'_{i,t}.$$
Then as an element in the graded piece $U_i^h/U_i^{h+1}$ (note that $T_i^{2h-1}dT_i = T_i^{2h}\dlog T_i \in V_i^{h+1}$ since $2h \ge h+1$ for $h\ge 1$, and we can pass from $V_i^{h+1}$ to $U_i^{h+1}$ by \eqref{modUvsmodV}),
\begin{align*} 
\eqref{prop:cl:eq12}
&=
\sum_{s\in I_i^{q-1}} 
T_i^{h}da'_{i,s}\wedge w_{i,s}
+
T_i^{h}\Big(
(-1)^{q-1}\sum_{s\in I_i^{q-1}} h a'_{i,s}w_{i,s}+
\sum_{t\in I_i^{q-2}}db'_{i,t}\wedge w'_{i,t}\Big)\wedge \dlog T_i
\\
&=\sum_{s\in I_i^{q-1}} 
T_i^{h}da_{i,s}\wedge \omega_{i,s}
+
T_i^{h}\Big(
(-1)^{q-1}\sum_{s\in I_i^{q-1}} h a_{i,s} \omega_{i,s}+
\sum_{t\in I_i^{q-2}}db_{i,t}\wedge \omega_{i,t}\Big)\wedge \dlog T_i
\\
&=
\rho_{i,h}\Big(
\sum_{s\in I_i^{q-1}}a_{i,s}\omega_{i,s},
\sum_{t\in I_i^{q-2}}b_{i,t}\omega_{i,t}
\Big)
\end{align*}
This proves the surjectivity.

An element
$(\sum a_{i,s}\omega_{i,s},
\sum b_{i,t}\omega_{i,t})$ lies in the kernel of \eqref{prop:cl:eq1} if and only if
\eq{prop:cl:eq2}{\sum_{s\in I_i^{q-1}} 
T_i^{h}
\cdot 
da_{i,s}\wedge\omega_{i,s}
+
T_i^{h}\cdot \Big(
(-1)^{q-1}\sum_{s\in I_i^{q-1}} ha_{i,s}\omega_{i,s}
+\sum_{t\in I_i^{q-2}}db_{i,t}\wedge\omega_{i,t}
\Big)\wedge 
\dlog T_i=0.}
If $p\nmid h$, this gives
\eq{lem3a1}{\Ker\eqref{prop:cl:eq1}
=
\Big\{\Big(
\frac{(-1)^{q}}{h}\sum_{t\in I_i^{q-2}}db_{i,t}\omega_{i,t},
\sum_{t\in I_i^{q-2}}b_{i,t}\omega_{i,t}
\Big)\mid 
b_{i,t}\in R_i
\Big\},
}
which is canonically isomorphic to $\Omega^{q-2}_{R_i}(\log A_i)$.
This proves the first statement in (1).
If $p\mid h$,
then \eqref{prop:cl:eq2} implies 
$d(\sum a_{i,s}\omega_{i,s})=0$, $d(\sum b_{i,t}\omega_{i,t})=0$, and therefore
$$\Ker\eqref{prop:cl:eq1}
=\cZ^{q-1} \oplus 
\cZ^{q-2} 
$$
This proves (2).

For the ``moreover'' part of (1), the surjectivity of the map $\alpha$ is clear.
\eqref{lem3a1} gives the exactness in the middle.
Suppose $\beta(\sum_{t}b_{i,t}\omega_{i,t})=0$ in $U_i^h/U_i^{h+1}$.
Note that
$$\beta(\sum_{t}b_{i,t}\omega_{i,t})
=\sum_t \dlog(1+b_{i,t}T_i^h)\wedge \omega_{i,t}\wedge\dlog T_i
\stackrel{\eqref{lem3a1}}{=}\sum_s \dlog(1+a_{i,s}T_i^h)\wedge\omega_{i,s}$$
for some $a_{i,s}$ with
$\sum_{s}a_{i,s}\omega_{i,s}
=-\frac{(-1)^{q}}{h}\sum_{t}db_{i,t}\omega_{i,t}$
in $\Omega^{q-1}_{R_i}(\log A_i)$.
Since we have proven that $\rho_{i,h}$ is an isomorphism, $\beta(\sum_{t}b_{i,t}\omega_{i,t})=0$ if and only if $\sum_s a_{i,s}\omega_{i,s}=0$. 
That is, $\sum_{t}b_{i,t}\omega_{i,t}\in \Ker d$.
This proves the remaining part of (1).
\item 
If $i\in [1,e]$ and $h=0$, any $q$-form in $U_i^{0}/U_i^{1}$ is of the form
$$
w+w'
\wedge \dlog T_i,
$$
where $w\in \Omega^q_{R_i}(\log A_i)_{\log},
w'\in \Omega^{q-1}_{R_i}(\log A_i)_{\log}$.
This proves (3).
\item 
Assume that $i\in [e+1,g]$ and $p\mid h$.
Based on \Cref{lem2}\eqref{lem2(2)}, a similar argument as in (a) shows that any $q$-form in $U_i^{h}/U_i^{h+1}$ is of the form
\begin{align*} 
\sum_{s\in I_i^{q-1}}
\dlog(1+a'_{i,s}T_i^{h+1})\wedge \omega_{i,s}
&+
\sum_{s\in I_i^{q-1}}
\dlog(1+a_{i,s}T_i^{h})\wedge \omega_{i,s}
\\
&+
\sum_{t\in I_i^{q-2}}
\dlog(1+b_{i,t}T_i^{h})\wedge \omega_{i,t}
\wedge \dlog (1+T_i),
\end{align*}
where $a'_{i,s},a_{i,s}, b_{i,t}\in R_i$.
Consider the map
\begin{align} 
\label{prop:cl:eq3}
\rho_{i,h}:
\Omega^{q-1}_{R_i}(\log A_i)
\oplus
\Omega^{q-1}_{R_i}(\log A_i)
\oplus
\Omega^{q-2}_{R_i}(\log A_i)
&\ra
U_i^{h}/U_i^{h+1}
\\
\Big(
\sum_{s\in I_i^{q-1}} a'_{i,s}\omega_{i,s},
\sum_{s\in I_i^{q-1}} a_{i,s}\omega_{i,s},
\sum_{t\in I_i^{q-2}} b_{i,t}\omega_{i,t}
\Big)
&\mapsto
\sum_{s\in I_i^{q-1}} \dlog(1+a'_{i,s}T_i^{h+1})\wedge\omega_{i,s}
\nonumber
\\
&\quad +
\sum_{s\in I_i^{q-1}} \dlog(1+a_{i,s}T_i^{h})\wedge\omega_{i,s}
\nonumber
\\
&\quad 
+\sum_{t\in I_i^{q-2}}
\dlog(1+b_{i,t}T_i^{h})\wedge\omega_{i,t}\wedge \dlog (1+T_i)
\nonumber
\end{align}
This map is then surjective.
An element
$(\sum a'_{i,s}\omega_{i,s},\sum a_{i,s}\omega_{i,s}
\sum b_{i,t}\omega_{i,t})$ lies in the kernel of \eqref{prop:cl:eq3} if and only if
$$\sum_{s\in I_i^{q-1}} 
T_i^{h}
\cdot 
da_{i,s}\wedge\omega_{i,s}
+ 
T_i^{h}\cdot \Big(
(-1)^{q-1}\sum_{s\in I_i^{q-1}} (h+1)a'_{i,s}\omega_{i,s}
+\sum_{t\in I_i^{q-2}}db_{i,t}\wedge\omega_{i,t}
\Big)\wedge 
\dlog (1+T_i)=0.
$$
Note that in this step, we have used the assumption that $h\ge 2$ and hence $a_{i,s}T_i^{h}d(a_{i,s}T_i^{h})\in V_i^{h+1}$.
As a result,
\begin{align}
\label{lem3c1}
\Ker\eqref{prop:cl:eq3}
&=
\Big\{\Big(
\frac{(-1)^{q}}{h+1}\sum_{t\in I_i^{q-2}}db_{i,t}\wedge\omega_{i,t},
\sum_{s\in I_i^{q-1}} a_{i,s}\omega_{i,s},
\sum_{t\in I_i^{q-2}}b_{i,t}\omega_{i,t}
\Big)
\mid 
d(\sum_{s\in I_i^{q-1}} a_{i,s}\omega_{i,s})=0
\Big\}
\end{align}
This proves the first half of (4).
The ``moreover'' part of (4) is proven analogously as the ``moreover'' part of (1).

\item 
Assume that $i\in [e+1,g]$ and $p\nmid h\ge p$.
Consider the map
\begin{align} 
\label{prop:cl:eq4}
\rho_{i,h}:
\Omega^{q-1}_{R_i}(\log A_i)
\oplus
\Omega^{q-2}_{R_i}(\log A_i)
&\ra
U_i^{h}/U_i^{h+1}
\\
\Big(
\sum_{s\in I_i^{q-1}} a_{i,s}\omega_{i,s},
\sum_{t\in I_i^{q-2}} b_{i,t}\omega_{i,t}
\Big)
&\mapsto
\sum_{s\in I_i^{q-1}} \dlog(1+a_{i,s}T_i^{h+1})\wedge\omega_{i,s}
\nonumber
\\
&\quad 
+\sum_{t\in I_i^{q-2}}
\dlog(1+b_{i,t}T_i^{h})\wedge\omega_{i,t}\wedge \dlog (1+T_i)
\nonumber
\end{align}
This map is surjective by a similar argument as in (a) based on \Cref{lem2}\eqref{lem2(3)}\eqref{lem2(4)}.
An element
$(\sum a_{i,s}\omega_{i,s},
\sum b_{i,t}\omega_{i,t})$ lies in the kernel of \eqref{prop:cl:eq4} if and only if
\eq{prop:cl:eq7}{
T_i^{h}\cdot \Big(
(-1)^{q-1} \sum_{s\in I_i^{q-1}} (h+1)a_{i,s}\omega_{i,s}
+\sum_{t\in I_i^{q-2}}db_{i,t}\wedge\omega_{i,t}
\Big)\wedge 
\dlog (1+T_i)=0.
}
If $p\nmid (h+1)$, \eqref{prop:cl:eq7} is equivalent to
$$\Ker\eqref{prop:cl:eq4}
=
\Big\{\Big(
\frac{(-1)^{q}}{h+1}
\sum_{t\in I_i^{q-2}}db_{i,t}\wedge\omega_{i,t},
\sum_{t\in I_i^{q-2}}b_{i,t}\omega_{i,t}
\Big)
\Big\}.
$$
This gives the first half of (5).
The ``moreover'' part of (5) is proven analogously as the ``moreover'' part of (1).
If $p\mid (h+1)$, 
\eqref{prop:cl:eq7} is equivalent to
$d(\sum b_{i,t}\omega_{i,t})=0,$
and hence
$$\Ker\eqref{prop:cl:eq4}
=\Omega^{q-1}_{R_i}(\log A_i)
\oplus 
\cZ^{q-2}
$$
This proves (6).

\item 
If $i\in [g+1,N]$ and $h=0$, any $q$-form in $U_i^{0}/U_i^{1}$ is of the form
\eq{prop:cl:eq5}{
w+w'
\wedge \dlog (1+T_i),
}
where $w\in \Omega^q_{R_i}(\log A)_{\log},
w'\in \Omega^{q-1}_{R_i}(\log A)_{\log}$.
This finishes (7).
\end{enumerate}
We hence finished the proof of \Cref{lem3}.
\end{proof}

\begin{lem}\label{lem4}
Let $w=\sum_{s\in I_i^{q-1}}a_{i,s}\omega_{i,s}\in \Omega_{R_i}^{q-1}(\log A_i)$,
$w'=\sum_{s\in I_i^{q-1}}a'_{i,s}\omega_{i,s}\in \Omega_{R_i}^{q-1}(\log A_i)$,
$w''=\sum_{s\in I_i^{q-2}}b_{i,t}\omega_{i,t}\in \Omega_{R_i}^{q-2}(\log A_i)$.
Let $h\ge 1$ be an integer. 
To simplify notations, we use
\begin{center} 
$\rho_{i,h}(w)$ (resp. $\rho_{i,h}(w,w'')$, $\rho_{i,h}(w',w)$ etc.) 
\end{center} 
to denote the elements in $U_i^{h}$ represented by the given form on right hand side of the isomorphisms in \Cref{lem3}.
\begin{enumerate}
\item 
\label{lem4(1)}
Suppose $i\in [1,e]$, $p\nmid h$ and $h\ge r_i$.    
If 
$\rho_{i,h}(w)\in G^{\ul r}_{R}\cap V_i^{h}\mod V_i^{h+1} ,$
then there exist
$$\tilde a_{i,s}, \tilde b_{i,t}\in (T_1^{\tilde r_1}\cdots \hat T_i^{\tilde r_i}\cdots T_N^{\tilde r_N})R_i$$
for all $s\in I_i^{q-1}$ and $t\in I_i^{q-2}$, such that
$$\rho_{i,h}(w)\equiv 
\sum_{s\in I_i^{q-1}} 
\dlog(1+\tilde a_{i,s}T_i^{h})\wedge \tilde \omega_{i,s}
+
\sum_{t\in I_i^{q-2}} \dlog(1+\tilde b_{i,t}T_i^{h})\wedge\tilde \omega_{i,t}\wedge \dlog T_i
\mod U_i^{h+1}.$$
\item \label{lem4(2)}
Suppose  $i\in [1,e]$, $p\mid h$ and $h\ge r_i$.    
If $\rho_{i,h}(w,w'')\in G^{\ul r}_{R}\cap V_i^{h}\mod V_i^{h+1},$
then there exist
$$\tilde a_{i,s}, \tilde b_{i,t}\in (T_1^{\tilde r_1}\cdots \hat T_i^{\tilde r_i}\cdots T_N^{\tilde r_N})R_i$$
such that 
$$\rho_{i,h}(w,w'')\equiv 
\sum_{s\in I_i^{q-1}} 
\dlog(1+\tilde a_{i,s}T_i^{h})\wedge \tilde \omega_{i,s}
+
\sum_{t\in I_i^{q-2}} 
\dlog(1+\tilde b_{i,t}T_i^{h})\wedge \tilde \omega_{i,t}\wedge\dlog T_i
\mod U_i^{h+1}.$$

\item \label{lem4(3)}
Suppose $i\in [e+1,g]$, $p\mid h$ and $h\ge r_i$.
If $\rho_{i,h}(w',w)\in G^{\ul r}_{R}\cap V_i^{h}\mod V_i^{h+2},$
then there exist
$$\tilde a'_{i,s}, \tilde a_{i,s}, \tilde b_{i,t}\in (T_1^{\tilde r_1}\cdots \hat T_i^{\tilde r_i}\cdots T_N^{\tilde r_N})R_i$$
for all $s\in I_i^{q-1}, t\in I_i^{q-2}$, such that
\begin{align*} 
\rho_{i,h}(w',w)
&\equiv 
\sum_{s\in I_i^{q-1}} \dlog(1+\tilde a'_{i,s}T_i^{h+1})\wedge\tilde\omega_{i,s}
+
\sum_{s\in I_i^{q-1}} \dlog(1+\tilde a_{i,s}T_i^{h})\wedge\tilde \omega_{i,s}
\\
&+ \sum_{t\in I_i^{q-2}} \dlog(1+\tilde b_{i,t}T_i^{h})\wedge\tilde \omega_{i,s}\wedge 
\dlog (1+T_i)
\mod U_i^{h+2}.
\end{align*}

\item \label{lem4(4)}
Suppose $i\in [e+1,g]$, $p\nmid h(h+1)$ and $h\ge r_i$.
If $\rho_{i,h}(w')\in G^{\ul r}_{R}\cap V_i^{h} \mod V_i^{h+2}$,
then there exist
$$
\tilde a'_{i,s},\tilde b_{i,t}\in (T_1^{\tilde r_1}\cdots \hat T_i^{\tilde r_i}\cdots T_N^{\tilde r_N})R_i$$
for all $s\in I_i^{q-1}, t\in I_i^{q-2}$, such that
$$\rho_{i,h}(w')
\equiv
\sum_{t\in I_i^{q-2}} 
\dlog(1+\tilde b_{i,t}T_i^{h})\wedge\tilde \omega_{i,t}\wedge \dlog (1+T_i)
+
\sum_{s\in I_i^{q-1}} \dlog(1+\tilde a'_{i,s}T_i^{h+1})\wedge\tilde\omega_{i,s}
\mod U_i^{h+1}
.$$

\item \label{lem4(5)}
Suppose $i\in [e+1,g]$,  $p\mid (h+1)$ and $h\ge r_i$.
If 
$\rho_{i,h}(w'')\in G^{\ul r}_{R}\cap V_i^{h} \mod V_i^{h+1}$,
then there exist
$$
\tilde b_{i,t}\in (T_1^{\tilde r_1}\cdots \hat T_i^{\tilde r_i}\cdots T_N^{\tilde r_N})R_i$$
for all $t\in I_i^{q-2}$, such that
$$\rho_{i,h}(w'')
\equiv 
\sum_{t\in I_i^{q-2}} 
\dlog(1+\tilde b_{i,t}T_i^{h})\wedge\tilde \omega_{i,t}
\wedge \dlog (1+T_i)
\mod U_i^{h+1}.$$

\end{enumerate}
\end{lem}

\begin{proof}[Proof of \Cref{lem4}]
(1)
Suppose $i\in [1,e]$, $p\nmid h$ and $h\ge r_i$.   
\Cref{lem3}\eqref{lem3(1} yields
\begin{align*} \rho_{i,h}(w)
&=\sum_{s\in I_i^{q-1}} \dlog(1+a_{i,s}T_i^{h})\wedge\omega_{i,s}
\\
&\equiv \sum_{s\in I_i^{q-1}}T_i^{h}
\left(
da_{i,s}\wedge \omega_{i,s}
    + (-1)^{q-1}  h a_{i,s} \omega_{i,s}\wedge \dlog T_i
\right)
\mod V_i^{h+1}.
\end{align*}
The condition $\rho_{i,h}(w)\in G_R^{\ul r} \cap V_i^{h}\mod V_i^{h+1}$ gives
$$\sum_{s\in I_i^{q-1}}T_i^{h}
\left(
da_{i,s}\wedge \omega_{i,s}
    + (-1)^{q-1}  h a_{i,s} \omega_{i,s}\wedge \dlog T_i
\right)\in G_R^{\ul r} \cap V_i^{h}$$
By \Cref{lem1}, 
\eq{lem4(1)-1}{
d(\sum_{s\in I_i^{q-1}} a_{i,s}\omega_{i,s})\in 
(T_1^{r_1}\cdots\hat T_i^{r_i}\cdots  T_N^{r_N})\Omega^q_{R_i}(\log A_i),
}
\eq{lem4(1)-2}{
\sum_{s\in I_i^{q-1}} a_{i,s}\omega_{i,s}\in 
(T_1^{r_1}\cdots\hat T_i^{r_i}\cdots  T_N^{r_N})\Omega^{q-1}_{R_i}(\log A_i).
}
By \Cref{RR24lem8.10gen}, the first relation implies that
there exist
$$\tilde a_{i,s}\in (T_1^{\tilde r_1}\cdots \hat T_i^{\tilde r_i}\cdots T_N^{\tilde r_N})R_i$$
for all $s$, such that
$d(\sum_{s} a_{i,s}\omega_{i,s})=d(\sum_{s} \tilde a_{i,s}\tilde \omega_{i,s})$.
The short exact sequence in \Cref{lem3}\eqref{lem3(1} implies 
$$\rho_{i,h}(w)-
\sum_{s\in I_i^{q-1}} \dlog(1+\tilde a_{i,s}T_i^{h})\wedge\tilde\omega_{i,s}
=\sum_{t\in I_i^{q-2}} \dlog(1+b_{i,t}T_i^{h})\wedge\omega_{i,s}\wedge \dlog T_i
$$
for some $b_{i,t}\in R_i$.
Since 
\begin{align*} 
\sum_{s\in I_i^{q-1}} \dlog(1+\tilde a_{i,s}T_i^{h})\wedge\tilde\omega_{i,s}
&=\sum_{s\in I_i^{q-1}} T_i^h
\left(d(\tilde a_{i,s}\tilde \omega_{i,s})
+(-1)^{q-1}\tilde a_{i,s}\tilde\omega_{i,s}\wedge \dlog T_i
\right)
\\
&=\sum_{s\in I_i^{q-1}} T_i^h
\left(d(a_{i,s} \omega_{i,s})
+(-1)^{q-1}\tilde a_{i,s}\tilde\omega_{i,s}\wedge\dlog T_i
\right)
\end{align*} 
and both $T_i^hd(a_{i,s} \omega_{i,s})$ and $T_i^h\tilde a_{i,s}\tilde\omega_{i,s}\wedge\dlog T_i$
lie in $G_R^{\ul r}\cap V_i^{h}$, we have
$$\sum_{s\in I_i^{q-1}} \dlog(1+\tilde a_{i,s}T_i^{h})\wedge\tilde\omega_{i,s}
\in G_R^{\ul r}\cap V_i^{h}.$$
From this and the assumption
$\rho_{i,h}(w)\in G_R^{\ul r}\cap V_i^{h} \mod V_i^{h+1}$, we deduce
$$\sum_{t\in I_i^{q-2}} \dlog(1+b_{i,t}T_i^{h})\wedge\omega_{i,t}\wedge \dlog T_i
\in G_R^{\ul r} \cap V_i^{h}
\mod V_i^{h+1}.$$
That is,
$$\sum_{t\in I_i^{q-2}} T_i^{h}db_{i,t}\wedge\omega_{i,t}\wedge \dlog T_i
\in G_R^{\ul r} \cap V_i^{h}$$
Again by \Cref{lem1}, it follows that
$$\sum_t d(b_{i,t}\omega_{i,t})\in (T_1^{r_1}\cdots\hat T_i^{r_i}\cdots  T_N^{r_N})\Omega^{q-1}_{R_i}(\log A_i).$$
\Cref{RR24lem8.10gen} gives the existence of 
$$\tilde b_{i,t}\in (T_1^{\tilde r_1}\cdots \hat T_i^{\tilde h}\cdots T_N^{\tilde r_N})R_i$$
for all $t$, such that
$d(\sum_t b_{i,t}\omega_{i,t})=d(\sum_t \tilde b_{i,t}\tilde \omega_{i,t}).$
Hence
$$\rho_{i,h}(w)\equiv
\sum_{s\in I_i^{q-1}} \dlog(1+\tilde a_{i,s}T_i^{h})\wedge\tilde\omega_{i,s}
+\sum_{t\in I_i^{q-2}} \dlog(1+\tilde b_{i,t}T_i^{h})\wedge\tilde \omega_{i,s}\wedge \dlog T_i
\mod V_i^{h+1}.
$$
This together with \eqref{modUvsmodV} proves (1).

\vspace{0.5em}
(2)
Suppose $i\in [1,e]$, $p\mid h$ and $h\ge r_i$.    
\Cref{lem3}\eqref{lem3(2} yields
\begin{align*} 
\rho_{i,h}(w,w'')
&=
\sum_{s\in I_i^{q-1}} 
\dlog(1+a_{i,s}T_i^{h})\wedge\omega_{i,s}
+
\sum_{t\in I_i^{q-2}} 
\dlog(1+b_{i,t}T_i^{h})\wedge\omega_{i,t}
\wedge \dlog T_i
\\
&\equiv 
\sum_{s\in I_i^{q-1}} 
T_i^{h}da_{i,s}
\wedge\omega_{i,s}
+
\sum_{t\in I_i^{q-2}} 
T_i^{h}db_{i,t}\wedge\omega_{i,t}
\wedge \dlog T_i
\mod V_i^{h+1}.
\end{align*}
The condition $\rho_{i,h}(w,w'')\in G_R^{\ul r} \cap V_i^{h}\mod V_i^{h+1}$ gives
$$\sum_{s\in I_i^{q-1}} 
T_i^{h}d(a_{i,s}
\omega_{i,s})
+
\sum_{t\in I_i^{q-2}} 
T_i^{h}d(b_{i,t}\omega_{i,t})
\wedge \dlog T_i
\in G_R^{\ul r} \cap V_i^{h}.$$
\Cref{lem1} implies that 
\eq{lem4(2)-1}{
d(\sum_{s\in I_i^{q-1}} a_{i,s}\omega_{i,s})\in 
(T_1^{r_1}\cdots\hat T_i^{r_i}\cdots  T_N^{r_N})\Omega^q_{R_i}(\log A_i),
}
\eq{lem4(2)-2}{
d(\sum_{t\in I_i^{q-2}} b_{i,t}\omega_{i,t})
\in (T_1^{r_1}\cdots\hat T_i^{r_i}\cdots T_N^{r_N})\Omega^{q-1}_{R_i}(\log A_i).
}
\Cref{RR24lem8.10gen} implies that there exist
$$\tilde a_{i,s}, \tilde b_{i,t}\in (T_1^{\tilde r_1}\cdots \hat T_i^{\tilde h}\cdots T_N^{\tilde r_N})R_i$$
such that 
$d(\sum_{s} a_{i,s}\omega_{i,s})=d(\sum_{s} \tilde a_{i,s}\tilde \omega_{i,s})$,
$d(\sum_{t} b_{i,t}\omega_{i,t})=d(\sum_{t} \tilde b_{i,t}\tilde \omega_{i,t})$.
Hence
\begin{align*}
\rho_{i,h}(w,w'')
&\equiv 
\sum_{s\in I_i^{q-1}} 
T_i^{h}d \tilde a_{i,s}
\wedge\tilde \omega_{i,s}
+
\sum_{t\in I_i^{q-2}} 
T_i^{h}d\tilde b_{i,t}\wedge\tilde \omega_{i,t}
\wedge \dlog T_i
\mod V_i^{h+1}
\\
&\equiv 
\sum_{s\in I_i^{q-1}} 
\dlog(1+\tilde a_{i,s}T_i^{h})\wedge \tilde \omega_{i,s}
+
\sum_{t\in I_i^{q-2}} 
\dlog(1+\tilde b_{i,t}T_i^{h})\wedge \tilde \omega_{i,t}\wedge\dlog T_i
\mod V_i^{h+1}.
\end{align*}
This together with \eqref{modUvsmodV} proves (2).

\vspace{0.5em}
(3)
Suppose $i\in [e+1,g]$, $p\mid h$ and $h\ge r_i$.
Note that these conditions imply that $h\ge 2$.
\Cref{lem3}\eqref{lem3(4} yields
\begin{align*}
\rho_{i,h}(
w',w
)
&= 
\sum_{s\in I_i^{q-1}} \dlog(1+a'_{i,s}T_i^{h+1})\wedge\omega_{i,s}
+
\sum_{s\in I_i^{q-1}} \dlog(1+a_{i,s}T_i^{h})\wedge\omega_{i,s}
\\
&\equiv
\sum_{s\in I_i^{q-1}} 
\left(
(-1)^{q-1}(h+1)
a'_{i,s}T_i^{h}\omega_{i,s}\wedge dT_i
+
T_i^{h+1}da'_{i,s}\wedge \omega_{i,s}
\right)
\\
&\quad +
\sum_{s\in I_i^{q-1}} T_i^{h} da_{i,s}\wedge\omega_{i,s}
\mod V_i^{h+2}
\\
&\equiv
(-1)^{q-1}(h+1)
\sum_{s\in I_i^{q-1}} 
a'_{i,s}T_i^{h}\omega_{i,s}\wedge \dlog (1+T_i)
+
\sum_{s\in I_i^{q-1}} T_i^{h} da_{i,s}\wedge\omega_{i,s}
\\
&\quad
+
(-1)^{q-1}(h+1)
\sum_{s\in I_i^{q-1}} 
a'_{i,s}T_i^{h+1}\omega_{i,s}\wedge \dlog (1+T_i)
+
\sum_{s\in I_i^{q-1}} 
T_i^{h+1}da'_{i,s}\wedge \omega_{i,s}
\mod V_i^{h+2}
.
\end{align*}
The condition $\rho_{i,h}(w,w'')\in G_R^{\ul r} \cap V_i^{h}\mod V_i^{h+2}$ gives
\begin{align*}
&(-1)^{q-1}(h+1)
\sum_{s\in I_i^{q-1}} 
a'_{i,s}T_i^{h}\omega_{i,s}\wedge \dlog (1+T_i)
+
\sum_{s\in I_i^{q-1}} T_i^{h} da_{i,s}\wedge\omega_{i,s}
\\
&
+
(-1)^{q-1}(h+1)
\sum_{s\in I_i^{q-1}} 
a'_{i,s}T_i^{h+1}\omega_{i,s}\wedge \dlog (1+T_i)
+
\sum_{s\in I_i^{q-1}} 
T_i^{h+1}da'_{i,s}\wedge \omega_{i,s}
\in G_R^{\ul r} \cap V_i^{h}.
\end{align*}
\Cref{lem1} implies that 
\eq{lem4(3)-1}{
\sum_{s\in I_i^{q-1}} a'_{i,s}\omega_{i,s}\in 
(T_1^{r_1}\cdots\hat T_i^{r_i}\cdots  T_N^{r_N})\Omega^{q-1}_{R_i}(\log A_i),
}
\eq{lem4(3)-2}{
d(\sum_{s\in I_i^{q-1}} a_{i,s}\omega_{i,s})\in 
(T_1^{r_1}\cdots\hat T_i^{r_i}\cdots  T_N^{r_N})\Omega^q_{R_i}(\log A_i).
}
\eq{lem4(3)-3}{
d(\sum_{s\in I_i^{q-1}} a'_{i,s}\omega_{i,s})\in 
(T_1^{r_1}\cdots\hat T_i^{r_i}\cdots  T_N^{r_N})\Omega^q_{R_i}(\log A_i).
}
By \Cref{RR24lem8.10gen}, \eqref{lem4(3)-2} implies
$d(\sum_{s} a_{i,s}\omega_{i,s})=d(\sum_{s} \tilde a_{i,s}\tilde \omega_{i,s})$ for some 
$$\tilde a_{i,s}\in (T_1^{\tilde r_1}\cdots\hat T_i^{\tilde h}\cdots T_N^{\tilde r_N}) R_i.$$
Similarly, \eqref{lem4(3)-3} implies
$d(\sum_{s} a'_{i,s}\omega_{i,s})=d(\sum_{s} \tilde a'_{i,s}\tilde \omega_{i,s})$ for some 
$$\tilde a'_{i,s}\in (T_1^{\tilde r_1}\cdots\hat T_i^{\tilde h}\cdots T_N^{\tilde r_N}) R_i.$$
In particular,
\begin{align*}
&
\rho_{i,h}(
w',w
)-
\sum_{s\in I_i^{q-1}} \dlog(1+\tilde a'_{i,s}T_i^{h+1})\wedge\tilde\omega_{i,s}
+
\sum_{s\in I_i^{q-1}} \dlog(1+\tilde a_{i,s}T_i^{h})\wedge\tilde \omega_{i,s}
\\
&\equiv 
(-1)^{q-1}(h+1)
\left(
\sum_{s\in I_i^{q-1}} 
(a'_{i,s}-\tilde a'_{i,s})T_i^{h}\tilde\omega_{i,s}
+
\sum_{s\in I_i^{q-1}} 
(a'_{i,s}-\tilde a'_{i,s})T_i^{h+1}\tilde\omega_{i,s}
\right)\wedge \dlog (1+T_i)
\\&\mod V_i^{h+2}.
\end{align*}
This is a closed form, so by the short exact sequence in \Cref{lem3}\eqref{lem3(4}, there exist some 
$b_{i,t}\in R_i$ for all $t$, such that
\begin{align*} 
&\rho_{i,h}(
w',w
)-
\sum_{s\in I_i^{q-1}} \dlog(1+\tilde a'_{i,s}T_i^{h+1})\wedge\tilde\omega_{i,s}
-
\sum_{s\in I_i^{q-1}} \dlog(1+\tilde a_{i,s}T_i^{h})\wedge\tilde \omega_{i,s}
\\
&\equiv \sum_{t\in I_i^{q-2}} \dlog(1+b_{i,t}T_i^{h})\wedge\omega_{i,s}\wedge 
\dlog (1+T_i)
\mod V_i^{h+2}.
\end{align*}
Since 
$\rho_{i,h}(
w',w
)$, $
\sum_{s\in I_i^{q-1}} \dlog(1+\tilde a'_{i,s}T_i^{h+1})\wedge\tilde\omega_{i,s}$ and 
$
\sum_{s\in I_i^{q-1}} \dlog(1+\tilde a_{i,s}T_i^{h})\wedge\tilde\omega_{i,s}$
lie in $G_R^{\ul r}\cap V_i^h \mod  V_i^{h+2}$,
we have
$$\sum_{t\in I_i^{q-2}} \dlog(1+b_{i,t}T_i^{h})\wedge\omega_{i,s}\wedge 
\dlog (1+T_i)
\in G_R^{\ul r}\cap V_i^h \mod  V_i^{h+2}$$
Taking $h\ge 2$ into consideration, a direct computation gives
\begin{align*}
\sum_{t\in I_i^{q-2}} \dlog(1+b_{i,t}T_i^{h})\wedge\omega_{i,t}\wedge \dlog (1+T_i)
&\equiv 
\sum_{t\in I_i^{q-2}}
T_i^h db_{i,t}\wedge\omega_{i,t}\wedge \dlog (1+T_i)
\mod V_i^{h+2}.
\end{align*}
Therefore
$$\sum_{t\in I_i^{q-2}}
T_i^h db_{i,t}\wedge\omega_{i,t}\wedge \dlog (1+T_i)
\in 
G_R^{\ul r}\cap V_i^h.
$$
Again by \Cref{lem1}, 
$$\sum_t d(b_{i,t}\omega_{i,t})\in (T_1^{r_1}\cdots\hat T_i^{r_i}\cdots  T_N^{r_N})\Omega^{q-1}_{R_i}(\log A_i).$$
\Cref{RR24lem8.10gen} gives the existence of 
$$\tilde b_{i,t}\in (T_1^{\tilde r_1}\cdots \hat T_i^{\tilde h}\cdots T_N^{\tilde r_N})R_i$$
for all $t$, such that
$d(\sum_t b_{i,t}\omega_{i,t})=d(\sum_t \tilde b_{i,t}\tilde \omega_{i,t}).$
In particular, 
$$\dlog(1+b_{i,t}T_i^{h})\wedge\omega_{i,t}\wedge \dlog (1+T_i)
\equiv 
\dlog(1+\tilde b_{i,t}T_i^{h})\wedge\tilde \omega_{i,t}\wedge \dlog (1+T_i)
\mod V_i^{h+2}.$$
Hence
\begin{align*}
\rho_{i,h}(w',w)
&\equiv
\sum_{s\in I_i^{q-1}} \dlog(1+\tilde a'_{i,s}T_i^{h+1})\wedge\tilde\omega_{i,s}
+
\sum_{s\in I_i^{q-1}} \dlog(1+\tilde a_{i,s}T_i^{h})\wedge\tilde\omega_{i,s}
\\
&+ \sum_{t\in I_i^{q-2}} \dlog(1+\tilde b_{i,t}T_i^{h})\wedge\tilde\omega_{i,s}\wedge 
\dlog (1+T_i)
\mod V_i^{h+2}.
\end{align*}
This together with \eqref{modUvsmodV} proves (3).

\vspace{0.5em}
(4)
Suppose $i\in [e+1,g]$, $p\nmid h(h+1)$ and $h\ge r_i$.
\Cref{lem3}\eqref{lem3(5} yields
\begin{align*} 
\rho_{i,h}(w')
&=
\sum_{s\in I_i^{q-1}} \dlog(1+a'_{i,s}T_i^{h+1})\wedge\omega_{i,s}\\
&\equiv
\sum_{s\in I_i^{q-1}} 
T_i^{h+1}da'_{i,s}\wedge \omega_{i,s}
+
(-1)^{q-1}(h+1)
\sum_{s\in I_i^{q-1}} 
a'_{i,s}T_i^{h}\omega_{i,s}\wedge dT_i
\mod V_i^{h+2}
\\
&=T_i^{h+1}
\left(
\sum_{s\in I_i^{q-1}} 
da'_{i,s}\wedge \omega_{i,s}
+
(-1)^{q-1}(h+1)
\sum_{s\in I_i^{q-1}} 
a'_{i,s}\omega_{i,s}\wedge \dlog (1+T_i)
\right)
\\
&+
(-1)^{q-1}(h+1)
\sum_{s\in I_i^{q-1}} 
a'_{i,s}T_i^{h}\omega_{i,s}\wedge \dlog (1+T_i)
\mod V_i^{h+2}.
\end{align*} 
In the equivalence relation step, we have used the relation $h+(h+1)\ge h+2$ as long as $h\ge 1.$
The condition
$\rho_{i,h}(w')
\in G_R^{\ul r}\cap V_i^{h}
\mod  V_i^{h+2}$ implies 
\begin{align*}
&T_i^{h+1}\left(
\sum_{s\in I_i^{q-1}} 
da'_{i,s}\wedge \omega_{i,s}
+
(-1)^{q-1}(h+1)
\sum_{s\in I_i^{q-1}} 
a'_{i,s}\omega_{i,s}\wedge \dlog (1+T_i)
\right)
\\
&+
(-1)^{q-1}(h+1)
\sum_{s\in I_i^{q-1}} 
a'_{i,s}T_i^{h}\omega_{i,s}\wedge \dlog (1+T_i)
\in G_R^{\ul r}\cap V_i^{h}.
\end{align*}
By \Cref{lem1}, this yields
\eq{lem4(4)-1}{\sum_{s\in I_i^{q-1}} 
a'_{i,s}\omega_{i,s}\wedge\dlog (1+T_i)\in 
(T_1^{r_1}\cdots\hat T_i^{r_i}\cdots  T_N^{r_N})\Omega^q_{R_i}(\log A_i).
}
\eq{lem4(4)-2}{\sum_{s\in I_i^{q-1}} 
a'_{i,s}\omega_{i,s}\wedge\dlog (1+T_i)\in 
(T_1^{r_1}\cdots\hat T_i^{r_i}\cdots  T_N^{r_N})\Omega^q_{R_i}(\log A_i).
}
\eq{lem4(4)-3}{
d\left(
\sum_{s\in I_i^{q-1}} 
a'_{i,s}\omega_{i,s}\right)
\in 
(T_1^{r_1}\cdots\hat T_i^{r_i}\cdots  T_N^{r_N})\Omega^q_{R_i}(\log A_i).
}
Applying \Cref{RR24lem8.10gen} to \eqref{lem4(4)-3}, we get some 
$$\tilde a'_{i,s} \in 
(T_1^{\tilde r_1}\cdots\hat T_i^{\tilde h}\cdots T_N^{\tilde r_N}) R_i$$
for all $s$, such that 
$d(\sum_{s}a'_{i,s}\omega_{i,s})=d(\sum_{s}\tilde a'_{i,s}\tilde \omega_{i,s})$.
The short exact sequence from \Cref{lem3}\eqref{lem3(5}, there exists some $b_{i,t}\in R_i$, such that
$$\rho_{i,h}(w')-
\sum_{s\in I_i^{q-1}} \dlog(1+\tilde a_{i,s}T_i^{h+1})\wedge\tilde\omega_{i,s}
=\sum_{t\in I_i^{q-2}} \dlog(1+b_{i,t}T_i^{h})\wedge\omega_{i,s}\wedge 
\dlog (1+T_i)
$$
Note that
$$\sum_{s\in I_i^{q-1}} \dlog(1+\tilde a_{i,s}T_i^{h+1})\wedge\tilde\omega_{i,s}
\in G_R^{\ul r}\cap V_i^{h}\mod V_i^{h+2}.$$
From this and the assumption
$\rho_{i,h}(w)\in G_R^{\ul r}\cap V_i^{h} \mod V_i^{h+2}$, we deduce
$$\sum_{t\in I_i^{q-2}} \dlog(1+b_{i,t}T_i^{h})\wedge\omega_{i,t}\wedge \dlog (1+T_i)
\in G_R^{\ul r} \cap V_i^{h}
\mod V_i^{h+2}.$$
A direct computation gives
\begin{align*}
\sum_{t\in I_i^{q-2}} \dlog(1+b_{i,t}T_i^{h})\wedge\omega_{i,t}\wedge \dlog (1+T_i)
\equiv \sum_{t\in I_i^{q-2}}
T_i^h db_{i,t}\wedge\omega_{i,t}\wedge \dlog (1+T_i)
\\
\begin{cases}
{\rm mod}\; V_i^{h+2},&h\ge 2\\
{\rm mod}\; V_i^{h+1},&h=1
\end{cases}
\end{align*}
Hence for any $h\ge 1$, the assumption
$\rho_{i,h}(w)\in G_R^{\ul r}\cap V_i^{h} \mod V_i^{h+2}$ forces 
$$\sum_{t\in I_i^{q-2}}
T_i^h db_{i,t}\wedge\omega_{i,t}\wedge \dlog (1+T_i)
\in G_R^{\ul r} \cap V_i^{h}.$$
By \Cref{lem1} again, 
$$\sum_t d(b_{i,t}\omega_{i,t})\in (T_1^{r_1}\cdots\hat T_i^{r_i}\cdots  T_N^{r_N})\Omega^{q-1}_{R_i}(\log A_i).$$
\Cref{RR24lem8.10gen} gives the existence of 
$$\tilde b_{i,t}\in (T_1^{\tilde r_1}\cdots \hat T_i^{\tilde h}\cdots T_N^{\tilde r_N})R_i$$
for all $t$, such that
$d(\sum_t b_{i,t}\omega_{i,t})=d(\sum_t \tilde b_{i,t}\tilde \omega_{i,t}).$
In particular, 
$$\dlog(1+b_{i,t}T_i^{h})\wedge\omega_{i,t}\wedge \dlog (1+T_i)
\equiv 
\dlog(1+\tilde b_{i,t}T_i^{h})\wedge\tilde\omega_{i,t}\wedge \dlog (1+T_i)
\;
\begin{cases}
{\rm mod}\; V_i^{h+2},&h\ge 2;\\
{\rm mod}\; V_i^{h+1},&h=1.
\end{cases}
$$
Hence
\begin{align*} 
\rho_{i,h}(w')
\equiv 
\sum_{t\in I_i^{q-2}} 
\dlog(1+\tilde b_{i,t}T_i^{h})\wedge\tilde\omega_{i,t}\wedge \dlog (1+T_i)
+
\sum_{s\in I_i^{q-1}} \dlog(1+\tilde a'_{i,s}T_i^{h+1})\wedge\tilde\omega_{i,s}
\\
\begin{cases}
{\rm mod}\; V_i^{h+2},&h\ge 2\\
{\rm mod}\; V_i^{h+1},&h=1
\end{cases}
\end{align*}
This together with \eqref{modUvsmodV} proves (a slightly stronger version of) \eqref{lem4(4)}.

\vspace{0.5em}
(5)
Suppose $i\in [e+1,g]$,  $p\mid (h+1)$ and $h\ge r_i$.
\begin{align*} 
\rho_{i,h}(w'')
&=
\sum_{t\in I_i^{q-2}} 
\dlog(1+b_{i,t}T_i^{h})\wedge\omega_{i,t}
\wedge \dlog (1+T_i)
\\
&\equiv
\sum_{t\in I_i^{q-2}} 
T_i^{h}db_{i,t}\wedge\omega_{i,t}
\wedge \dlog (1+T_i)
\mod V_i^{h+1},
\end{align*}
The condition
$\rho_{i,h}(w')\in G_R^{\ul r}\cap V_i^{h}\mod V_i^{h+1}$ implies 
$$\sum_{t\in I_i^{q-2}} 
T_i^{h}d(b_{i,t}\omega_{i,t})
\wedge \dlog (1+T_i)
\in G_R^{\ul r}\cap V_i^{h}.$$
\Cref{lem1} gives
\eq{lem4(5)-1}{
d(\sum_{t\in I_i^{q-2}} b_{i,t}\omega_{i,t})
\in (T_1^{r_1}\cdots\hat T_i^{r_i}\cdots  T_N^{r_N})\Omega^{q-1}_{R_i}(\log A_i).}
\Cref{RR24lem8.10gen} implies that $d(\sum_{t\in I_i^{q-2}} b_{i,t}\omega_{i,t})=d(\sum_{t\in I_i^{q-2}} \tilde b_{i,t}\tilde \omega_{i,t})$ for some 
$$\tilde b_{i,t}\in 
(T_1^{\tilde r_1}\cdots\hat T_i^{\tilde h}\cdots T_N^{\tilde r_N}) R_i.
$$
Moreover,
\begin{align*} 
\rho_{i,h}(w'')&
\equiv
\sum_{t\in I_i^{q-2}} 
T_i^{h}d(b_{i,t}\omega_{i,t})
\wedge \dlog (1+T_i)
\mod V_i^{h+1}
\\
&\equiv
\sum_{t\in I_i^{q-2}} 
T_i^{h}d(\tilde b_{i,t}\tilde \omega_{i,t})
\wedge \dlog (1+T_i)
\mod V_i^{h+1}
\\
&\equiv 
\sum_{t\in I_i^{q-2}} 
\dlog(1+\tilde b_{i,t}T_i^{h})\wedge\tilde \omega_{i,t}
\wedge \dlog (1+T_i)
\mod V_i^{h+1}.
\end{align*}
This together with \eqref{modUvsmodV} proves (5).
\end{proof}

We continue with the proof of \Cref{MainProp}.
For any $i\in [1,g]$, any $w\in G^{\ul r}_{R,\log}$  can be written as an $T_i$-adically convergent infinite sum 
\eq{MainPropconteq}{w_{i,r_i}+w_{i,r_i+1}+w_{i,r_i+2}+\cdots}
in the module of continuous K\"ahler differentials $\lim_n \Omega^q_{(R/(T_1,\cdots ,T_N)^n)/k'}$,
where $w_{i,l}$ (for each $l\ge r_i$) is the \textit{homogeneous degree $l$ part} of $w$ in the $i$th direction in our convention, namely,
the projection of $w$ to
$$ 
(T_1^{r_1}\cdots T_i^{l}\cdots T_N^{r_N})\cdot
\Omega^q_{R_i}(\log A_i)\,\oplus \,
(T_1^{r_1}\cdots T_i^{l}\cdots T_N^{r_N})\cdot
\Omega^{q-1}_{R_i}(\log A_i)\dlog T_i,
$$
in case $i\in [1,e]$, and to
$$ 
(T_1^{r_1}\cdots T_i^{l}\cdots T_N^{r_N})\cdot
\Omega^q_{R_i}(\log A_i)\,\oplus \,
(T_1^{r_1}\cdots T_i^{l}\cdots T_N^{r_N})\cdot
\Omega^{q-1}_{R_i}(\log A_i)\dlog (1+T_i)
$$
in case $i\in [e+1,g]$.
In particular,
$$w_{i,l}\in  G_R^{\ul r}\cap U_i^{l}.$$
Since $k$ is perfect, the continuous K\"ahler differentials canonically agrees with the algebraic ones. 
Therefore $w$ agrees with the infinite sum \eqref{MainPropconteq}.

If $e>0$, pick an $i\in [1,e]$. 
We do an induction definition of $w_{i,l}^{\rm st}$ on $l$.
Set $w_{i,0}^{\rm st}=\cdots=w_{i,r_i-1}^{\rm st}=0$.
Now suppose that $w_{i,r_i-1}^{\rm st},\cdots, w_{i,l-1}^{\rm st}$ are defined such that
$$w\equiv \sum_{j=r_i}^{l-1} w_{i,j}^{\rm st}\mod U_i^{l-2}.$$
We define an homogeneous degree $l$ element $w_{i,l}^{\rm st}$ in $U_i^l$.
By \Cref{lem4}\eqref{lem4(1)}(2), 
there exist
$$a_{i,s,l}\in (T_1^{\tilde r_1}\cdots \hat T_i^{\tilde r_i}\cdots T_N^{\tilde r_N})\cdot R_i,\quad 
\text{for every }s\in I_i^{q-1}\text{ and every integer } l\ge r_i$$
and 
$$b_{i,t, l}\in (T_1^{\tilde r_1}\cdots \hat T_i^{\tilde r_i}\cdots T_N^{\tilde r_N})\cdot R_i,\quad
\text{for every }t\in I_i^{q-2}\text{ and every integer } l\ge r_i,$$
such that
$$w - \sum_{j=r_i}^{l-1} w_{i,j}^{\rm st}
\equiv
w_{i,l}^{\rm st}
\mod U_i^{l+1}$$
where 
$$w_{i,l}^{\rm st}
=\sum_{s\in I_i^{q-1}}\sum_{l\ge r_i}
\dlog(1+a_{i,s,l}T_i^{l})\wedge \tilde \omega_{i,s}
+
\sum_{t\in I_i^{q-2}} \sum_{l\ge r_i}
\dlog(1+b_{i,t,l}T_i^{l})\wedge \tilde \omega_{i,t}\wedge \dlog T_i
.$$
Altogether,
\begin{align*}
w&=
\sum_{s\in I_i^{q-1}}\sum_{l\ge r_i}
\dlog(1+a_{i,s,l}T_i^{l})\wedge \tilde \omega_{i,s}
+
\sum_{t\in I_i^{q-2}}\sum_{l\ge r_i}
\dlog (1+b_{i,t,l}T_i^{l})
\wedge \tilde \omega_{i,t}\wedge \dlog T_i
\\
&=
\sum_{s\in I_i^{q-1}}\dlog x\wedge \tilde \omega_{i,s}
+
\sum_{t\in I_i^{q-2}}
\dlog y
\wedge \tilde \omega_{i,t}\wedge \dlog T_i
\end{align*}
with
$$x=\prod_{l\ge r_i} (1+a_{i,s,l}T_i^l),
\quad
y=\prod_{l\ge r_i} (1+b_{i,t,l}T_i^{l}).
$$
(When $q=1$, we have $y=1$ as a convention.)
Then both $x$ and $y$ converge in the $T_i$-adic, and therefore the $(T_1,\dots,T_N)$-adic topology of $R$. 
Moreover, 
$$x\in 1+(T_1^{\tilde r_1}\cdots T_N^{\tilde r_N}),
\quad 
y\in 1+(T_1^{\tilde r_1}\cdots  T_N^{\tilde r_N}).$$
This completes the proof of \Cref{MainProp} in the ``$\subset$'' direction in this case.

If $e=0$, we choose an $i\in [e+1,g]$.
We do an induction definition of $w_{i,l}^{\rm st}$ on $l$.
Set $w_{i,0}^{\rm st}=\cdots=w_{i,r_i-1}^{\rm st}=0$.
Now suppose that $w_{i,r_i-1}^{\rm st},\cdots, w_{i,l-1}^{\rm st}$ are defined such that
$$w\equiv \sum_{j=r_i}^{l-1} w_{i,j}^{\rm st}\mod U_i^{l-2}.$$
We define an homogeneous degree $l$ element $w_{i,l}^{\rm st}$ in $U_i^l$.
\begin{itemize}
\item 
Case 1. $l=pm$ for some $m\ge 1$.

By \Cref{lem3}\eqref{lem3(4},
there exist $a'_{i,s,l}$, $a_{i,s,l}\in R_i$ for all $s$, such that
\begin{align*}
w'_{i,l}:=&\sum_{s\in I_i^{q-1}} \dlog(1+a'_{i,s,l}T_i^{l+1})\wedge\omega_{i,s}
+
\sum_{s\in I_i^{q-1}} \dlog(1+a_{i,s,l}T_i^{l})\wedge\omega_{i,s}
\end{align*} 
represent the same residue class as $w- \sum_{j=r_i}^{l-1} w_{i,j}^{\rm st}$ in $U_i^{l}/U_i^{l+1}$.
By \Cref{lem4}\eqref{lem4(3)}, 
there exist $\tilde a'_{i,s,l}$, $\tilde a_{i,s,l}$, $\tilde b_{i,t,l}\in 
(T_1^{\tilde r_1}\cdots \hat T_i^{\tilde r_i}\cdots T_N^{\tilde r_N})\cdot R_i$ for all $s,t$ such that
\begin{align*}
w_{i,l}^{\rm st}:=
&\sum_{s\in I_i^{q-1}} 
\dlog(1+\tilde a'_{i,s,l}T_i^{l+1})\wedge\tilde\omega_{i,s}
+
\sum_{s\in I_i^{q-1}} 
\dlog(1+\tilde a_{i,s,l}T_i^{l})\wedge\tilde \omega_{i,s}
\\
&+ \sum_{t\in I_i^{q-2}} 
\dlog(1+\tilde b_{i,t,l}T_i^{l})\wedge\tilde\omega_{i,s}\wedge 
\dlog (1+T_i)
\end{align*}
represents the same residue class as $w'_{i,l}$ in $U_i^{l}/U_i^{l+2}$.

\item 
Case 2. $p\nmid l(l+1)$.

By \Cref{lem3}\eqref{lem3(5},
there exist $a'_{i,s,l}\in R_i$ for all $s$ such that
\begin{align*}
w'_{i,l}:=&\sum_{s\in I_i^{q-1}} \dlog(1+a'_{i,s,l}T_i^{l+1})\wedge\omega_{i,s}
\end{align*} 
represent the same residue class as $w- \sum_{j=r_i}^{l-1} w_{i,j}^{\rm st}$ in $U_i^{l}/U_i^{l+1}$.
By \Cref{lem4}\eqref{lem4(4)}, 
there exist $\tilde a'_{i,s,l}$,  $\tilde b_{i,t,l}\in
(T_1^{\tilde r_1}\cdots \hat T_i^{\tilde r_i}\cdots T_N^{\tilde r_N})\cdot R_i$  for all $s,t$ such that
\begin{align*}
w_{i,l}^{\rm st}:=
&\sum_{s\in I_i^{q-1}} 
\dlog(1+\tilde a'_{i,s,l}T_i^{l+1})\wedge\tilde\omega_{i,s}
+ \sum_{t\in I_i^{q-2}} 
\dlog(1+\tilde b_{i,t,l}T_i^{l})\wedge\tilde\omega_{i,s}\wedge 
\dlog (1+T_i)
\end{align*}
represents the same residue class as $w'_{i,l}$ in $U_i^{l}/U_i^{l+1}$.

\item 
Case 3. $p\mid l+1.$

By \Cref{lem3}\eqref{lem3(6},
there exist $b_{i,t,l}\in R_i$ for all $t$ such that
\begin{align*}
w'_{i,l}:=\sum_{t\in I_i^{q-2}} 
\dlog(1+ b_{i,t,l}T_i^{l})\wedge\omega_{i,s}\wedge 
\dlog (1+T_i)
\end{align*} 
represent the same residue class as $w- \sum_{j=r_i}^{l-1} w_{i,j}^{\rm st}$ in $U_i^{l}/U_i^{l+1}$.
By \Cref{lem4}\eqref{lem4(5)}, 
there exist $\tilde b_{i,t,l}\in 
(T_1^{\tilde r_1}\cdots \hat T_i^{\tilde r_i}\cdots T_N^{\tilde r_N})\cdot R_i$ for all $t$, such that
\begin{align*}
w_{i,l}^{\rm st}:=
\sum_{t\in I_i^{q-2}} 
\dlog(1+\tilde b_{i,t,l}T_i^{l})\wedge\tilde \omega_{i,s}\wedge 
\dlog (1+T_i)
\end{align*}
represents the same residue class as $w'_{i,l}$ in $U_i^{l}/U_i^{l+1}$.
\end{itemize}
Since in every case, $w'_{i,l}$ is the homogeneous degree $l$ part of $w- \sum_{j=r_i}^{l-1} w_{i,j}^{\rm st}$,
the difference
$$w- \sum_{j=r_i}^{l} w_{i,j}^{\rm st}$$
lies in $U_i^{l+1}$. This finishes the induction step.

Now consider the sum
$$w^{\rm st}:=\sum_{l\ge r_i} w_{i,l}^{\rm st}. $$
This is a convergent differential form in the module of continuous Kähler differentials.
In fact, set 
$$l_0=\begin{cases}
    r_i,&\text{if $p\mid r_i$;}
    \\
    r_i+1,&\text{if $p\nmid r_i$.}
\end{cases}$$
To sum it up: there exist
$$a_{i,s,l}\in (T_1^{\tilde r_1}\cdots \hat T_i^{\tilde r_i}\cdots T_N^{\tilde r_N})\cdot R_i,\quad 
\text{for every }s\in I_i^{q-1}\text{ and every  }  l\ge l_0,$$
$$b_{i,t, l}\in (T_1^{\tilde r_1}\cdots \hat T_i^{\tilde r_i}\cdots T_N^{\tilde r_N})\cdot R_i,\quad
\text{for every }t\in I_i^{q-2}\text{ and every $l\ge r_i$,}$$
such that
\begin{align*} 
w&=
\sum_{s\in I_i^{q-1}}\sum_{l\ge l_0} \dlog(1+a_{i,s,l}T_i^{l})\wedge\tilde \omega_{i,s}
+
\sum_{t\in I_i^{q-2}}\sum_{l\ge r_i}
\dlog(1+ b_{i,t,l}T_i^{l})\wedge\tilde \omega_{i,t}\wedge \dlog (1+ T_i)
\\
&=
\sum_{s\in I_i^{q-1}}\dlog x\wedge\tilde \omega_{i,s}
+
\sum_{t\in I_i^{q-2}}
\dlog y
\wedge \tilde \omega_{i,t}\wedge \dlog (1+T_i)
\end{align*}
with
$$x=\prod_{l\ge l_0} (1+a_{i,s,l}T_i^{l}),
\quad 
y=\prod_{l\ge r_i} (1+b_{i,t,l}T_i^{l}).
$$
It is therefore clear that
$$x\in 1+(T_1^{\tilde r_1}\cdots T_N^{\tilde r_N}),
\quad 
y\in 1+(T_1^{\tilde r_1}\cdots  T_N^{\tilde r_N}).$$
This finishes the proof of \Cref{MainProp} in the ``$\subset$'' direction.

\section{Appendix}
\subsection{A Hartog's type lemma for Cohen-Macaulay abelian sheaves}
\begin{lem}
\label{lemAppHartog}
Let $X$ be a locally noetherian equidimensional scheme, and
$F$ be a Cohen-Macaulay abelian sheaf on $X_\Zar$ (see \cite[Definition and Remark on p.238]{Ha66}), and $U$ be an open dense subset of $X$. Then for any section $\alpha\in \Gamma(U,F)$ such that $\alpha\in F_x$ for all $x\in X^{(1)}$, we have $\alpha\in \Gamma(X,F).$ 

\end{lem}
\begin{proof}
If $I^\bullet$ is an injective resolution of $F$, then 
$$0\ra \ul\Gamma_{\bar{\{x\}}}(I^\bullet) \ra I^\bullet\ra j_*I^\bullet
\ra 0,$$
is a short exact of complexes, where $j:X\setminus \bar{\{x\}}\hra X$ is the open immersion. 
The long exact sequence it induces gives a four-term exact sequence of sheaves
$$0\ra \ul\Gamma_{\bar{\{x\}}}(F) \ra F \ra j_*F\xra{\delta} \ul H^1_{\bar{\{x\}}}(F)\ra 0.$$
Taking stalks at $x$, we get a canonical isomorphism
\eq{lemAppHartog1}{H^1_{x}(F)\simeq \bigoplus_{\stackrel{\eta\in X^{(0)}}{x\in\bar{\{\eta\}}}} F_{\eta}/F_{x}}
for every $x\in X^{(1)}$.
Since $F$ is Cohen-Macauley, the Cousin complex
$$0\ra 
\bigoplus_{\eta\in X^{(0)}} i_{\eta*}F_\eta \ra
\bigoplus_{x\in X^{(1)}} i_{x*}  H_x^1(F)\ra\cdots
$$
is a flasque resolution of the Zariski sheaf $F$.
Bringing \eqref{lemAppHartog1} into this sequence, we deduce that
$$0\ra 
F\ra 
\bigoplus_{\eta\in X^{(0)}} i_{\eta*}F_\eta \ra
\bigoplus_{x\in X^{(1)}} i_{x*}\Big(\bigoplus_{\stackrel{\eta\in X^{(0)}}{x\in\bar{\{\eta\}}}} F_{\eta}/F_{x}\Big)
$$
is exact.
Applying the functor $\Gamma(X,-)$, we get a left exact sequence
$$0\ra 
\Gamma(X,F)\ra 
\bigoplus_{\eta\in X^{(0)}} \Gamma(X_\eta ,F_\eta) \ra
\bigoplus_{x\in X^{(1)}} \Big(\bigoplus_{\stackrel{\eta\in X^{(0)}}{x\in\bar{\{\eta\}}}} F_{\eta}/F_{x}\Big).
$$
Hence if $\alpha\in F_x$ for all $x\in X^{(1)}$, then the left exactness of the above sequence gives $\alpha\in \Gamma(X,F)$.
\end{proof}
\begin{example}
\label{lemAppHartog:Ex}
Let $k$ be a perfect field of positive characteristic. The followings are examples of Cohen-Macaulay sheaves:
\begin{enumerate}
\item 
\label{lemAppHartog:Ex1}
Locally free sheaves of finite rank on a smooth $k$-scheme are Cohen-Macaulay. (\cite[p.239, Example]{Ha66})
\item 
\label{lemAppHartog:Ex2}
For any $n\ge 1$ and $q\ge 0$, $W_n\Omega^q_X$ is a Cohen-Macaulay $W_n\cO_X$-module when $X$ is smooth over $k$ by\cite[I. Corollaire 3.9]{IlDRW}. It is not a locally free $W_n\cO_X$-module.
\item 
\label{lemAppHartog:Ex3}
More generally, $W_n\Omega^q_{(X,\pm D)}$ as defined in \cite{RamFil1} is a Cohen-Macaulay $W_n\cO_X$-module when $X$ is separated smooth of finite type over $k$.  
This follows from \cite[Theorem 6.4]{RamFil1} and the injectivity of $\ul p$ for the pair $(X,D)$, and from the duality statement Theorem 9.3 in \textit{loc. cit.} for the pair $(X,-D)$.
\end{enumerate}
\end{example}

\begin{cor}\label{lemAppHartogNis}
Let $X$ be a locally noetherian integral equidimensional normal scheme, and $F$ is a locally free sheaf of $\cO_X$-modules of finite rank on $X_\Nis$, whose restriction to $X_\Zar$ is Cohen-Macaulay.
Then for any section $\alpha\in \Gamma(U,F)$ such that $\alpha\in F_{x,\Nis}$ (= the Nisnevich stalk) for all $x\in X^{(1)}$, we have $\alpha\in \Gamma(X,F).$
\end{cor}
\begin{proof}
For any $x\in X^{(1)}$, the Nisnevich stalk of the structure sheaf is a henselian discrete valuation ring by \cite[\href{https://stacks.math.columbia.edu/tag/0AP3}{Tag 0AP3}]{stacks-project}.
In particular, we have
$$\cO_{x,\Zar}=\cO_{x,\Nis}\cap \Frac \cO_{X,x}.$$
As $F$ is a locally free $\cO_X$-module of finite rank, say $r$, we have
$$F_{x,\Zar}=F_{x,\Nis}\cap (\Frac \cO_{X,x})^r.$$
As $F|_{X_\Zar}$ is a Cohen-Macaulay Zariski sheaf,  the Corollary follows from \Cref{lemAppHartog}.
\end{proof}

\subsection{{A generalization of \cite[Lemma 8.10]{RamFil1}}}
\begin{lem}
\label{RR24lem8.10gen}
Let $X$ be a separated smooth $k$-scheme of finite type with \'etale coordinates $T_1,\cdots, T_d\in \Gamma(X,\cO_X)$.
Set 
$$A=\Div(T_1\cdots T_e), \quad
B=\Div(T_1^{r_1}\cdots   T_g^{r_g}),$$
$$\tilde A=\Div(T_1\cdots T_e T_{e+1}\cdots T_{f}),
\quad 
\tilde B=\Div(T_1^{\tilde r_1}\cdots  T_d^{\tilde r_d}),$$
where $\ul r'=(\tilde r_1, \cdots, \tilde r_d)$ is defined by
\eq{rprimedef}{
\tilde r_j=r_j \text{ for all $j\notin [e+1,f]$, 
\quad and \quad
$\tilde r_j=r_j+1$ for all $j\in [e+1,f]$.}
}
Then
$dF_{X*}\left(\Omega^{q-1}_{X}(\log \tilde A)(-\tilde B)\right)$ is a locally free sheaf of finite rank, and
we have an identification of Nisnevich sheaves 
\begin{align*} 
\left( \Omega^q_{X}(\log A)(-B)\right)\cap d\Omega^{q-1}_X(\log A)
=
d\left(\Omega^{q-1}_{X}(\log \tilde A)(-\tilde B)\right)
\end{align*}
on $X_\Nis$.
\end{lem}
\begin{proof}[Proof of \Cref{RR24lem8.10gen}]
We first show the local freeness of $dF_{X*}(\Omega^{q-1}_X(\log \tilde A)(-\tilde B))$.
The definition of $\tilde A$ and $\tilde B$ guarantees that any $T_j$ from $\tilde B$ that does not lie in $\tilde A$ has its power divisible by $p$.
Therefore to show $dF_{X*}(\Omega^{q-1}_X(\log \tilde A)(-\tilde B))$ is locally free, we can assume $f=g$.
With this extra assumption, the desired local freeness follows from an decreasing argument on $q$ via the following two exact sequences, where $\cZ^q$ is defined by the exactness of the first:
$$0\ra \cZ^{q-1}\ra 
\Omega_X^{q-1}(\log \tilde A)(-\tilde B)\xra{d} 
d\Omega_X^{q-1}(\log \tilde A)(-\tilde B)
\ra 0 $$
$$0\ra 
d\Omega_X^{q-1}(\log A)(-\tilde B)\ra 
\cZ^q \xra{C} 
\Omega_X^{q}(\log \tilde A)(-\lceil \frac{\tilde B}{p}\rceil ) 
\ra 0.$$

The second part of this lemma is a generalization of \cite[Lemma 8.10]{RamFil1}.
The direction ``$\supset$'' is clear.
Now we show ``$\subset$''.
Because $dF_{X*}(\Omega^{q-1}_X(\log \tilde A)(-\tilde B))$ is locally free,
it suffices to show the inclusion for the Nisnevich stalk at any codimension 1 point $x$ by \Cref{lemAppHartogNis}.
The statement is only nontrivial when $x$ is the generic point of the reduced divisors cut out by $T_j$ for $j\in [1,g]$.
If $j\in [1,e]$, or $j\in [f+1,g]$, the lemma follows from the part (a) and part (b) of the proof of \cite[Lemma 8.10]{RamFil1}.
Now we suppose $x$ is the generic point of $\Div (T_j)$ for some $j\in [e+1,f]$.
Write $T:=T_j$ and $r:=r_j$. 
By replacing $X$ with a Nisnevich neighborhood $\Spec R$ of $\Div(T)$, we can assume $R_0:=R/T$ is isomorphic to a subring of $R$ (\cite[Lemma 7.14]{BRS}).
Because we are looking at the Nisnevich stalk at the generic point of $\Div(T)$, we can assume that $X$ is of the form $\Spec R_0[T]$.
Hence as free $R_0$-modules, 
$$
\Omega^{q-1}_X(\log A)=\Omega^{q-1}_{R_0[T]}=
\bigoplus_{i\ge 0} \Omega^{q-1}_{R_0}\cdot T^i
\oplus 
\bigoplus_{i\ge 1} \Omega^{q-2}_{R_0} \cdot T^i\dlog T,
$$
$$\Omega^{q-1}_X(\log \tilde A)=\Omega^{q-1}_{R_0[T]}(\log T)=
\bigoplus_{i\ge 0} \Omega^{q-1}_{R_0}\cdot T^i
\oplus 
\bigoplus_{i\ge 0} \Omega^{q-2}_{R_0} \cdot T^i\dlog T.$$
Now suppose $\alpha\in \Omega^{q-1}_{X}(\log A)$.
Then there exist
$a_i\in \Omega^{q-1}_{R_0}, b_i\in \Omega^{q-2}_{R_0}$, among which only finitely many are nonzero, such that
$$\alpha=a_0+ \sum_{i\ge 1} (a_iT^i+ b_iT^i\dlog T).$$
Hence
$$d\alpha=\sum_{i\ge 0}T^ida_i+\sum_{i\ge 1}((-1)^{q-1}ia_{i}+ db_i)\cdot T^{i}\dlog T.$$
The condition $d\alpha\in T^r \Omega^q_{R_0[T]}
=\bigoplus_{i\ge r} \Omega^q_{R_0}\cdot T^{i}\oplus \bigoplus_{i\ge r+1} \Omega^{q-1}_{R_0} \cdot T^{i}\dlog T$ 
gives
$$da_i=0\text{  for all $i\in [0,r-1],$
\quad and 
\quad 
$(-1)^{q-1}ia_{i}+db_i=0$ for all $i\in [1,r]$.}$$
In particular, the second equation above implies $db_i=0$ when $p\mid i\in [1,r]$; and $da_r=0$ (since $p\nmid r$ by our assumption).
Let 
$$\epsilon:= a_0+\sum_{p\mid i\in [1,r]} a_iT^i+ \sum_{i\in [1,r]} b_iT^i\dlog T
-\sum_{p\nmid i\in [1,r]} \frac{(-1)^{q-1}}{i} T^idb_i.$$
It is direct to check that 
$$d\epsilon=da_0+\sum_{p\mid i\in [1,r]} T^i\left(da_i+db_i\dlog T\right)=0.$$
Moreover,
$$\alpha-\epsilon
=
\sum_{i\ge r+1}(a_iT^i+b_iT^i \dlog T)
\in T^{r+1}\cdot \Omega^{q-1}_X(\log \tilde A).$$
\end{proof}

\subsection{Log forms as the Frobenius invariant in $W_n\Omega^q_X(\log A)/\sB_\infty$}
The following lemma originates from a lemma by Kay Rülling. 
Below we formulate and prove a log version.
We only need the $n=1$ version of it, but nevertheless we record a
proof for general $n$.

\begin{lem}
\label{sesF-1Binf}
Let $X$ be a connected smooth $k$-scheme and $A$ be a simple normal crossing divisor on $X$. 
Let $j:U:=X\setminus A\hra X$ be the open immersion.
Let $W_n\Omega^*_{X}(\log A)$ be the de Rham-Witt complex associated to the log structure on $X$ associated to $A$ (relative to the trivial log structure on the base $k$) defined by Hyodo-Kato.
More explicitly (see \cite[Proposition 4.6]{HK}),
$$W_n\Omega^q_{X}(\log A)
=
W_n\cO_X\cdot j_*\dlog K^M_{q,U}+dW_n\cO_X\cdot j_*\dlog K^M_{q-1,U}
\subset j_*W_n\Omega^q_{U},$$
where $K^M_{q,U}$ is the \'etale sheaf of Milnor $K$-theory on $U$.
Let 
$$W_n\Omega^q_{X}(\log A)_{\log}:=W_n\Omega^q_{X}(\log A)\cap j_*W_n\Omega^q_{U,\log}.$$
Then the sequence 
$$0\ra W_n\Omega^q_{X}(\log A)_{\log}\ra W_n\Omega^q_X(\log A)/\sB_\infty\xra{\bar F-1} W_n\Omega^q_X(\log A)/\sB_\infty\ra 0$$
is exact on $X_\et.$
Here 
\footnote{Note that by \cite[I.3.11.4]{IlDRW}, for $n=1$, the image of $F^j d W_{j+1}\Omega^{q-1}$ in $\Omega^q$ is exactly $\sB_{j+1}\Omega^q$, so this definition recovers $\sB_\infty = \bigcup_{m\ge 1} \sB_m\Omega^q$ as used previously.}
$$\sB_\infty:=\bigcup_{j\ge 0}F^j(dW_{n+j}\Omega^{q-1}_U(\log A)),$$
and $\bar F$ is induced from $F$ via the restriction map.
\end{lem}
\begin{proof}
We first show that the map $\bar F$ is well-defined.
By \cite[Theorem 4.4]{HK}, 
$$\Ker(R:W_{n+1}\Omega^q_X(\log A)\ra W_n\Omega^q_X(\log A))
=V^n(\Omega^q_X(\log A))+dV^n(\Omega^{q-1}_X(\log A)).
$$
Hence $F:W_{n+1}\Omega^q_U\ra W_n\Omega^q_U$ induces a well-defined map
$$\bar F: W_n\Omega^q_X(\log A)\ra W_n\Omega^q_X(\log A)/dV^{n-1}\Omega^{q-1}_X(\log A).$$
Since $dV^{n-1}\Omega^{q-1}_X(\log A)\subset dW_n\Omega^{q-1}_X(\log A)\subset \sB_\infty$, the map $F$ further induces a well-defined map
$$W_n\Omega^q_X(\log A)/\sB_\infty\ra W_n\Omega^q_X(\log A)/\sB_\infty$$
which we also denote by $\bar F$.

Clearly the sequence is a complex.
When $A=\emptyset$, the surjectivity of $\bar F-1$  follows from \cite[I, Proposition 3.26]{IlDRW}.
In general, we can run a similar (but way easier) argument as in the proof in \textit{loc. cit.}. 
For any $\alpha\in W_n\cO_X\cdot j_*\dlog K^M_{q-1,U}$, 
we have $\sum_{i=1}^n dV^i(R^{i-1}\alpha)
\in W_{n+1}\Omega^q_X(\log A)$
and it satisfies
$$d\alpha=(\bar F-1)(R\Big(\sum_{i=1}^n dV^i(R^{i-1}\alpha)\Big)).$$
It remains to show $W_n\cO_X\cdot j_*\dlog K^M_{q,U}\subset (\bar F-1)(W_n\Omega^q_X(\log A))$ \'etale locally.
Since $\bar F-1:W_n\cO_X\ra W_n\cO_X$ is surjective on $X_\et$ and log forms are invariant under the endomorphism $\bar F$, 
it follows that $W_n\cO_X\cdot j_*\dlog K^M_{q,U}\subset (\bar F-1)(W_n\cO_X\cdot j_*\dlog K^M_{q,U})$, and in particular , it lies in $ (\bar F-1)(W_n\Omega^q_X(\log A))$ as well.


We continue to prove the injectivity.
Let $\alpha \in W_{n}\Omega^q_{X}(\log A)_{\log}$ be zero in $W_n\Omega^q_X(\log A)/\sB_\infty$. 
Then by the definition of $\sB_\infty$,
$\alpha=F^m d\beta$ for some $\beta\in W_{n+m}\Omega^{q-1}_X(\log A)$.
It follows that 
\begin{align*}
V^{n+m}(\alpha) 
&= V^{n+m}(F^m d\beta) 
=p^mV^n d\beta
= p^{n+m} d V^n \beta 
= dV^n(p^{n+m}\beta)=0.
\end{align*}
Let $\tilde \alpha\in W_{2n+m}\Omega^q_X(\log A)_{\log}$ be a lift of $\alpha$.
Then 
$$0= V^{n+m}(\alpha)= V^{n+m}(F^{n+m}(\tilde \alpha))= \ul p^{n+m}(\alpha),$$
where in the second equality we have used that $\tilde \alpha$ is a log form and hence $F(\tilde \alpha)=R(\tilde \alpha)$.
By the injectivity of $\ul p$, we deduce that $\alpha=0$ as desired.

Next we prove the exactness in the middle.
Let $\alpha\in W_{n}\Omega^q_{X}(\log A)$  and assume  
$\bar F(\alpha)-\alpha = F^m d\beta$
for some $\beta\in W_{n+m}\Omega_X^{q-1}(\log A)$.
Let 
$$\gamma=\sum_{i=1}^{n+m} dV^i(R^{i-1}(\beta))\;\in W_{n+m+1}\Omega^q_X(\log A).$$
Then 
$d\beta=(\bar F-1)(R\gamma)$, and
$$\bar F(\alpha-F^m R(\gamma)) - (\alpha- F^m R(\gamma))=
F^md\beta -\bar F(F^{m}R\gamma)+F^m R(\gamma)
=0.$$
By \cite[Lemme 2]{CTSS}, $\alpha- F^m R(\gamma)$ is a log form.
\end{proof}
\section*{Conflict of interest statement }
On behalf of all authors, the corresponding author states that there is no conflict of interest.

\bibliographystyle{alpha}
\bibliography{Cond-DRW}

\end{document}